\documentclass{elsart}
\usepackage{natbib}
\bibpunct{(}{)}{;}{a}{,}{,}
\usepackage{amssymb}
\usepackage{amsmath}
\usepackage{pstricks}
\usepackage{pst-node}

\newcommand{\fso}{\mathfrak{so}}
\newcommand{\fsp}{\mathfrak{sp}}
\newcommand{\fsl}{\mathfrak{sl}}
\newcommand{\bZ}{\mathbb{Z}}
\newcommand{\bQ}{\mathbb{Q}}
\newcommand{\cA}{\mathcal{A}}
\DeclareMathOperator{\Hom}{Hom}
\DeclareMathOperator{\End}{End}
\DeclareMathOperator{\supp}{supp}
\newcommand{\sr}[1]{\stackrel{#1}{\longrightarrow}}

\SpecialCoor
\newcommand{\singleloop}{\pspicture[.4](-.6,-.5)(.6,.5)
\pscircle(0,0){.4}
\endpspicture}
\newcommand{\doubleloop}{\pspicture[.4](-.6,-.5)(.6,.5)
\pscircle[doubleline=true](0,0){.4}
\endpspicture}
\newcommand{\pentanode}{\pnode(.4;90){a1}\pnode(.4;162){a2}\pnode(.4;234){a3}
\pnode(.4;306){a4}\pnode(.4;18){a5}
\pnode(.9;90){b1}\pnode(.9;162){b2}\pnode(.9;234){b3}
\pnode(.9;306){b4}\pnode(.9;18){b5}}
\newcommand{\hexanode}{\pnode(.4;0){a1}\pnode(.4;60){a2}\pnode(.4;120){a3}
\pnode(.4;180){a4}\pnode(.4;240){a5}\pnode(.4;300){a6}
\pnode(.9;0){b1}\pnode(.9;60){b2}\pnode(.9;120){b3}
\pnode(.9;180){b4}\pnode(.9;240){b5}\pnode(.9;300){b6}}

\journal{Advances in Mathematics}
\begin{document}
\begin{frontmatter}
\title{Invariant tensors for the spin representation of $\mathfrak{so}(7)$}
\author{Bruce W. Westbury}
\address{Mathematics Institute\\
University of Warwick\\
Coventry CV4 7AL}
\ead{bww@maths.warwick.ac.uk}
\date{\today}
\begin{abstract} We construct a pivotal category by a finite
presentation and show that it is an integral form of the
category of invariant tensors of the spin representation of
the quantum group $U_q(B_3)$ over the field $\bQ(q)$.
\end{abstract}
\begin{keyword}invariant tensors, pivotal category
\end{keyword}
\end{frontmatter}

\section{Introduction}
In this paper we describe the category of invariant tensors of the
rank three simple Lie algebra $\fso(7)=B_3$ as a finitely presented
pivotal category; and this paper can be considered as a response to
\cite[Question 3.6]{MR1265145}. This follows the work of Greg Kuperberg in
\cite{MR1403861} which describes the categories of invariant tensors
for the rank two simple Lie algebras as finitely presented pivotal
categories. The most interesting case is the exceptional Lie algebra
$G_2$ which is discussed in \cite{MR1265145} and in
\cite{math.CO/0507112}. This work was preceeded by the Temperley-Lieb
category which gives the category of invariant tensors for the 
rank one simple Lie algebra $\fsl(2)$: a full account of this is
given in \cite{MR1446615}.

The spin representation of $\mathfrak{so}(7)$ is discussed in \cite{MR1999634}
and \cite{MR1048745} and $R$-matrices are discussed in \cite{MR1121816} and
\cite{MR1086739}. The main aim of this paper to to give a definition
of a finitely presented strict pivotal category by using diagrams
and the main result is that this category is isomorphic as a strict
pivotal category to the category of invariant tensors for this spin
representation.

This representation has the property that all non-zero weight spaces are
one dimensional. This implies that the tensor product with any highest
weight representation is multiplicity free. In particular the tensor
product with itself is multiplicity free. This tensor product is the
sum of four representations. These representations are the trivial
representation, the vector representation, the adjoint representation
and the third exterior power of the vector representation.

The (quantum) dimensions are
\begin{center}
\begin{tabular}{cccc}
Vector & Adjoint & Spin & \\
$[1,0,0]$ & $[0,1,0]$ & $[0,0,1]$ & $[0,0,2]$ \\ \hline
& & & \\
$\frac{[10][7]}{[5][2]}$ & $\frac{[12][7][6]}{[4][3][2]}$ &%
$\frac{[10][6][2]}{[5][3][1]}$  & $\frac{[12][10][7]}{[6][4][1]}$ \\
\end{tabular}
\end{center}

The organisation of this paper is that in \S \ref{diagrams} we give
the diagrammatic relations which are relations for a finitely 
presented $\bZ[\delta]$-linear strict pivotal category. 
In \S \ref{fpres} we define a sequence of finitely presented
algebras. In \S \ref{crystal} we give a diagrammatic model for
the crystal graphs of the spin representation, the vector
representation and their tensor products. In \S \ref{conc}
we relate these constructions and identify the sequence of
finitely presented algebras with the centraliser algebras of the
tensor powers of the spin representation. In \S \ref{conf} we
prove that the presentations of the centraliser algebras are
confluent and in \S \ref{cell} we show that these centraliser
algebras are cellular.

\section{Diagrams}\label{diagrams}
In this section we give a finite presentation of a $\bZ[\delta]$-linear
strict pivotal category.
There are two types of edge which we draw as a single line and
a double line. There is just one type of vertex, namely the
trivalent vertex with two single lines and one double line.
$$
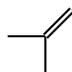
\pspicture(-.7,-.7)(.7,.7)
\qline(-.5,0)(0,0)\qline(0,-.5)(0,0)\psline[doubleline=true](0,0)(.35,.35)
\endpspicture
$$

The coefficients in the following equations are written in
terms of the quantum integers. Although they may be written as
a fraction these rational functions are all Laurent polynomials
in $q$ invariant under the involution $q\leftrightarrow q^{-1}$.
This ring can be identified with the polynomial ring $\bZ[\delta]$
by the inclusion $\delta\mapsto q+q^{-1}$. The following tables can be
used to convert all the coefficients that appear to elements of
$\bZ[\delta]$.
\begin{center}
\begin{tabular}{cccc}
$[1]$ & $[2]$ & $[3]$ & $[4]$ \\ \hline
$1$ & $\delta$ & $\delta^2-1$ & $\delta(\delta^2-2)$ \\
\end{tabular}
\begin{tabular}{cccc}
$\frac{[4]}{[2]}$ & $\frac{[6]}{[3]}$ &$\frac{[8]}{[4]}$ &$\frac{[10]}{[5]}$ \\
\hline
$\delta^2-2$ & $\delta(\delta^2-3)$ & $\delta^4-4\delta^2+2$ &%
$\delta(\delta^4-5\delta^2+5)$ \\
\end{tabular}
\end{center}
These two sequences of polynomials, $\{P(k)\}$, each satisfies the 
recurrence relation $P(k+1)-\delta P(k)+P(k-1)=0$
(but with different initial conditions).

Define a diagram to be a trivalent graph embedded in a rectangle with
boundary points on the top and bottom edges.
Each edge is labelled as being either single or double such that
there are two single edges and one double edge at each trivalent
vertex. Let $\widetilde{\mathcal{D}}$ be the category whose morphisms are
isotopy equivalence classes of such diagrams. The composition of morphisms
is given by putting one rectangle on top of the other.
This category is a strict pivotal category. The tensor product
of morphisms is given by putting the rectangles together side by
side and the dual diagram is obtained by rotating through half a
revolution.

This category is referred to as the free strict pivotal category
on the trivalent vertex since it has the following universal
property. Let $\mathcal{P}$ be a strict pivotal category with
self-dual objects $V$ and $W$ and a morphism 
$V\otimes V\rightarrow W$. Then there is a unique functor of
strict pivotal categories $\widetilde{\mathcal{D}}\rightarrow\mathcal{P}$
which maps the trivalent vertex to this morphism.

Then we take the free $\bZ[\delta]$-linear category on
$\widetilde{\mathcal{D}}$ and impose relations.
This defines the category $\mathcal{D}^\prime$. This means that
$\mathcal{D}^\prime$
is a $\bZ[\delta]$-linear strict pivotal category and that there is
a functor of $\bZ[\delta]$-linear strict pivotal categories
$\widetilde{\mathcal{D}}\rightarrow\mathcal{D}^\prime$.
The defining relations consist of 
the basic relations given in Figure \ref{relb},
the two square relations given in Figure \ref{rels},
the pentagon relation given in Figure \ref{relp}
and the hexagon relation given in Figure \ref{relh}.

The objects of $\widetilde{\mathcal{D}}$ and $\mathcal{D}^\prime$
are finite sequences of single and double edges. The category
$\mathcal{D}$ is the full subcategory of $\mathcal{D}^\prime$
whose objects are finite sequences of single edges.
\begin{figure}
\begin{eqnarray}
\singleloop & = & \frac{[10][6][2]}{[5][3]} \nonumber \\
\doubleloop & = & \frac{[10][7]}{[5][2]} \nonumber \\
\pspicture[.4](-.6,-.5)(.6,.5)
\psline[doubleline=true](-.5,0)(0,0)
\psbezier(0,0)(.7,.7)(.7,-.7)(0,0)
\endpspicture
& = & 0 \nonumber \\
\pspicture[.4](-.8,-.5)(.8,.5)
\psline[doubleline=true](-.7,0)(-.3,0)
\pcarc[arcangle=45](-.3,0)(.3,0)
\pcarc[arcangle=-45](-.3,0)(.3,0)
\psline[doubleline=true](.3,0)(.7,0)
\endpspicture
& = & \frac{[2]^2[6]}{[3]}\pspicture[.4](-.6,-.5)(.6,.5)
\psline[doubleline=true](-.5,0)(.5,0)
\endpspicture \nonumber \\
\pspicture[.4](-.8,-.5)(.8,.5)
\psline(-.7,0)(-.3,0)
\pcarc[arcangle=45,doubleline=true](-.3,0)(.3,0)
\pcarc[arcangle=-45](-.3,0)(.3,0)
\psline(.3,0)(.7,0)
\endpspicture
& = & [7]\pspicture[.4](-.6,-.5)(.6,.5)
\psline(-.5,0)(.5,0)
\endpspicture \nonumber \\
\pspicture[.4](-.9,-.9)(.9,.9)
\psline[doubleline=true](.4;90)(.8;90)
\psline[doubleline=true](.4;210)(.8;210)
\psline[doubleline=true](.4;330)(.8;330)
\pcarc[arcangle=-15](.4;90)(.4;210)
\pcarc[arcangle=-15](.4;210)(.4;330)
\pcarc[arcangle=-15](.4;330)(.4;90)
\endpspicture
& = & 0 \nonumber \\
\pspicture[.4](-.9,-.9)(.9,.9)
\psline[doubleline=true](.4;90)(.8;90)
\psline(.4;210)(.8;210)
\psline(.4;330)(.8;330)
\pcarc[arcangle=-15](.4;90)(.4;210)
\pcarc[arcangle=-15,doubleline=true](.4;210)(.4;330)
\pcarc[arcangle=-15](.4;330)(.4;90)
\endpspicture
& = & -[5]
\pspicture[.4](-.9,-.9)(.9,.9)
\psline[doubleline=true](.0;90)(.8;90)
\psline(.0;210)(.8;210)
\psline(.0;330)(.8;330)
\endpspicture\nonumber 
\end{eqnarray}
\caption{Basic relations}\label{relb}
\end{figure}

\begin{figure}
\begin{eqnarray}
\pspicture[.4](-.9,-.9)(.9,.9)
\psline(.4;45)(.8;45)\psline(.4;135)(.8;135)
\psline(.4;225)(.8;225)\psline(.4;315)(.8;315)
\pcarc[arcangle=-15,doubleline=true](.4;45)(.4;135)
\pcarc[arcangle=-15](.4;135)(.4;225)
\pcarc[arcangle=-15,doubleline=true](.4;225)(.4;315)
\pcarc[arcangle=-15](.4;315)(.4;45)
\endpspicture
& = & 
[3] \pspicture[.4](-.6,-.5)(.6,.5)
\psbezier(.5;45)(.25;45)(.25;135)(.5;135)
\psbezier(.5;225)(.25;225)(.25;315)(.5;315)
\endpspicture
+[4] \pspicture[.4](-.6,-.5)(.6,.5)
\psbezier(.5;45)(.25;45)(.25;315)(.5;315)
\psbezier(.5;225)(.25;225)(.25;135)(.5;135)
\endpspicture
+\frac{[4]}{[2]} \pspicture[.4](-.8,-.5)(.8,.5)
\qline( .25,0)( .5,.433)\qline( .25,0)( .5,-.433)
\qline(-.25,0)(-.5,.433)\qline(-.25,0)(-.5,-.433)
\psline[doubleline=true](-.25,0)(.25,0)
\endpspicture
+[4] \pspicture[.4](-.6,-.7)(.6,.7)
\qline(0, .25)(.433, .5)\qline(0, .25)(-.433, .5)
\qline(0,-.25)(.433,-.5)\qline(0,-.25)(-.433,-.5)
\psline[doubleline=true](0,-.25)(0,.25)
\endpspicture \label{sq1} \\
\pspicture[.4](-.9,-.9)(.9,.9)
\psline[doubleline=true](.4;45)(.8;45)
\psline[doubleline=true](.4;135)(.8;135)
\psline(.4;225)(.8;225)\psline(.4;315)(.8;315)
\pcarc[arcangle=-15](.4;45)(.4;135)\pcarc[arcangle=-15](.4;135)(.4;225)
\pcarc[arcangle=-15,doubleline=true](.4;225)(.4;315)
\pcarc[arcangle=-15](.4;315)(.4;45)
\endpspicture
& = & [2]^2 \pspicture[.4](-.6,-.5)(.6,.5)
\psbezier[doubleline=true](.5;45)(.25;45)(.25;135)(.5;135)
\psbezier(.5;225)(.25;225)(.25;315)(.5;315)
\endpspicture
+ [3] \pspicture[.4](-.8,-.5)(.8,.5)
\psline[doubleline=true]( .25,0)( .5,.433)\qline( .25,0)( .5,-.433)
\psline[doubleline=true](-.25,0)(-.5,.433)\qline(-.25,0)(-.5,-.433)
\qline(-.25,0)(.25,0)
\endpspicture
\label{sq2}
\end{eqnarray}
\caption{Square relations}\label{rels}
\end{figure}

\begin{figure}
\begin{multline}
\pspicture[.4](-1,-1)(1,1)\pentanode
\ncline[doubleline=true]{a1}{b1}\ncline{a2}{b2}\ncline{a3}{b3}\ncline{a4}{b4}
\ncline{a5}{b5}
\ncarc[arcangle=6]{a1}{a5}
\ncarc[arcangle=6]{a2}{a1}
\ncarc[arcangle=6,doubleline=true]{a3}{a2}
\ncarc[arcangle=6]{a4}{a3}
\ncarc[arcangle=6,doubleline=true]{a5}{a4}
\endpspicture
 =  -\left(
\pspicture[.4](-1,-1)(1,1)\pentanode
\ncline[doubleline=true]{a1}{b1}\ncline{a2}{b2}\ncline{a3}{b3}
\ncarc[arcangle=6,doubleline=true]{a3}{a2}
\ncarc[arcangle=6]{a2}{a1}
\nccurve[angleA= 24,angleB=126]{a3}{b4}
\nccurve[angleA=330,angleB=198]{a1}{b5}
\endpspicture +
\pspicture[.4](-1,-1)(1,1)\pentanode\ncline[doubleline=true]{b1}{a1}
\nccurve[angleA=210,angleB=342]{a1}{b2}
\nccurve[angleA= 54,angleB=126]{b3}{b4}
\nccurve[angleA=330,angleB=198]{a1}{b5}
\endpspicture +
\pspicture[.4](-1,-1)(1,1)\pentanode
\ncline{a4}{b4}\ncline{a5}{b5}\ncline[doubleline=true]{a1}{b1}
\ncarc[arcangle=6]{a1}{a5}
\ncarc[arcangle=6,doubleline=true]{a5}{a4}
\nccurve[angleA=240,angleB=342]{a1}{b2}
\nccurve[angleA=186,angleB= 54]{a4}{b3}
\endpspicture\right) \\
 -[2] \left(
\pspicture[.4](-1,-1)(1,1)\pentanode
\ncline{a2}{b2}\ncline{a3}{b3}\ncline{a4}{b4}
\ncarc[arcangle=6,doubleline=true]{a4}{a3}
\ncarc[arcangle=6]{a3}{a2}
\nccurve[angleA= 96,angleB=198]{a4}{b5}
\nccurve[angleA= 42,angleB=270,doubleline=true]{a2}{b1}
\endpspicture +
\pspicture[.4](-1,-1)(1,1)\pentanode
\ncline{a3}{b3}\ncline{a4}{b4}\ncline{a5}{b5}
\ncarc[arcangle=6]{a5}{a4}
\ncarc[arcangle=6,doubleline=true]{a4}{a3}
\nccurve[angleA=168,angleB=270,doubleline=true]{a5}{b1}
\nccurve[angleA=114,angleB=342]{a3}{b2}
\endpspicture +
\pspicture[.4](-1,-1)(1,1)\pentanode\ncline{b2}{a2}
\nccurve[angleA=282,angleB= 54]{a2}{b3}
\nccurve[angleA=126,angleB=198]{b4}{b5}
\nccurve[angleA= 42,angleB=270,doubleline=true]{a2}{b1}
\endpspicture +
\pspicture[.4](-1,-1)(1,1)\pentanode\ncline{b5}{a5}
\nccurve[angleA=138,angleB=270,doubleline=true]{a5}{b1}
\nccurve[angleA=342,angleB= 54]{b2}{b3}
\nccurve[angleA=258,angleB=126]{a5}{b4}
\endpspicture
\right) \label{pt}
\end{multline}
\caption{Pentagon relation}\label{relp}
\end{figure}

\begin{figure}
\begin{multline}
\pspicture[.4](-1,-1)(1,1)\hexanode
\ncline{a1}{b1}\ncline{a2}{b2}\ncline{a3}{b3}
\ncline{a4}{b4}\ncline{a5}{b5}\ncline{a6}{b6}
\ncarc[arcangle=6]{a1}{a6}
\ncarc[arcangle=6,doubleline=true]{a2}{a1}
\ncarc[arcangle=6]{a3}{a2}
\ncarc[arcangle=6,doubleline=true]{a4}{a3}
\ncarc[arcangle=6]{a5}{a4}
\ncarc[arcangle=6,doubleline=true]{a6}{a5}
\endpspicture 
+ \left\{
\pspicture[.4](-1,-1)(1,1)\hexanode
\ncline{a1}{b1}\ncline{a2}{b2}\ncline{a3}{b3}
\ncline{a4}{b4}\ncline{a5}{b5}\ncline{a6}{b6}
\ncarc[arcangle=6]{a2}{a1}
\ncarc[arcangle=6,doubleline=true]{a3}{a2}
\ncarc[arcangle=6]{a4}{a3}
\ncarc[arcangle=6,doubleline=true]{a5}{a4}
\ncarc[arcangle=6]{a6}{a5}
\endpspicture
+
\pspicture[.4](-1,-1)(1,1)\hexanode
\ncline{a1}{b1}\ncline{a2}{b2}\ncline{a3}{b3}
\ncline{a4}{b4}\ncline{a5}{b5}\ncline{a6}{b6}
\ncarc[arcangle=6]{a1}{a6}
\ncarc[arcangle=6,doubleline=true]{a2}{a1}
\ncarc[arcangle=6]{a3}{a2}
\ncarc[arcangle=6]{a5}{a4}
\ncarc[arcangle=6,doubleline=true]{a6}{a5}
\endpspicture
+
\pspicture[.4](-1,-1)(1,1)\hexanode
\ncline{a1}{b1}\ncline{a2}{b2}\ncline{a3}{b3}
\ncline{a4}{b4}\ncline{a5}{b5}\ncline{a6}{b6}
\ncarc[arcangle=6]{a1}{a6}
\ncarc[arcangle=6,doubleline=true]{a2}{a1}
\ncarc[arcangle=6]{a3}{a2}
\ncarc[arcangle=6,doubleline=true]{a4}{a3}
\ncarc[arcangle=6]{a5}{a4}
\endpspicture \right \} \\
+ \left\{
\pspicture[.4](-1,-1)(1,1)\hexanode
\ncline{a1}{b1}
\ncline{a4}{b4}\ncline{a5}{b5}\ncline{a6}{b6}
\ncarc[arcangle=6]{a1}{a6}
\ncarc[arcangle=6]{a5}{a4}
\ncarc[arcangle=6,doubleline=true]{a6}{a5}
\nccurve[angleA=240,angleB=300]{b2}{b3}
\endpspicture
+
\pspicture[.4](-1,-1)(1,1)\hexanode
\ncline{a1}{b1}\ncline{a2}{b2}\ncline{a3}{b3}
\ncline{a6}{b6}
\ncarc[arcangle=6]{a1}{a6}
\ncarc[arcangle=6,doubleline=true]{a2}{a1}
\ncarc[arcangle=6]{a3}{a2}
\nccurve[angleA=0,angleB=60]{b4}{b5}
\endpspicture
+
\pspicture[.4](-1,-1)(1,1)\hexanode
\ncline{a2}{b2}\ncline{a3}{b3}
\ncline{a4}{b4}\ncline{a5}{b5}
\ncarc[arcangle=6]{a3}{a2}
\ncarc[arcangle=6]{a5}{a4}
\ncarc[arcangle=6,doubleline=true]{a4}{a3}
\nccurve[angleA=180,angleB=120]{b1}{b6}
\endpspicture 
\right\} +
\pspicture[.4](-1,-1)(1,1)\hexanode
\nccurve[angleA=240,angleB=300]{b2}{b3}
\nccurve[angleA=0,angleB=60]{b4}{b5}
\nccurve[angleA=120,angleB=180]{b6}{b1}
\endpspicture \\
=
\pspicture[.4](-1,-1)(1,1)\hexanode
\ncline{a1}{b1}\ncline{a2}{b2}\ncline{a3}{b3}
\ncline{a4}{b4}\ncline{a5}{b5}\ncline{a6}{b6}
\ncarc[arcangle=6,doubleline=true]{a1}{a6}
\ncarc[arcangle=6]{a2}{a1}
\ncarc[arcangle=6,doubleline=true]{a3}{a2}
\ncarc[arcangle=6]{a4}{a3}
\ncarc[arcangle=6,doubleline=true]{a5}{a4}
\ncarc[arcangle=6]{a6}{a5}
\endpspicture
+ \left\{
\pspicture[.4](-1,-1)(1,1)\hexanode
\ncline{a1}{b1}\ncline{a2}{b2}\ncline{a3}{b3}
\ncline{a4}{b4}\ncline{a5}{b5}\ncline{a6}{b6}
\ncarc[arcangle=6]{a2}{a1}
\ncarc[arcangle=6,doubleline=true]{a3}{a2}
\ncarc[arcangle=6]{a4}{a3}
\ncarc[arcangle=6,doubleline=true]{a5}{a4}
\ncarc[arcangle=6]{a6}{a5}
\endpspicture +
\pspicture[.4](-1,-1)(1,1)\hexanode
\ncline{a1}{b1}\ncline{a2}{b2}\ncline{a3}{b3}
\ncline{a4}{b4}\ncline{a5}{b5}\ncline{a6}{b6}
\ncarc[arcangle=6,doubleline=true]{a1}{a6}
\ncarc[arcangle=6]{a2}{a1}
\ncarc[arcangle=6]{a4}{a3}
\ncarc[arcangle=6,doubleline=true]{a5}{a4}
\ncarc[arcangle=6]{a6}{a5}
\endpspicture +
\pspicture[.4](-1,-1)(1,1)\hexanode
\ncline{a1}{b1}\ncline{a2}{b2}\ncline{a3}{b3}
\ncline{a4}{b4}\ncline{a5}{b5}\ncline{a6}{b6}
\ncarc[arcangle=6,doubleline=true]{a1}{a6}
\ncarc[arcangle=6]{a2}{a1}
\ncarc[arcangle=6,doubleline=true]{a3}{a2}
\ncarc[arcangle=6]{a4}{a3}
\ncarc[arcangle=6]{a6}{a5}
\endpspicture \right\} \\
+ \left\{
\pspicture[.4](-1,-1)(1,1)\hexanode
\ncline{a3}{b3}
\ncline{a4}{b4}\ncline{a5}{b5}\ncline{a6}{b6}
\ncarc[arcangle=6]{a4}{a3}
\ncarc[arcangle=6,doubleline=true]{a5}{a4}
\ncarc[arcangle=6]{a6}{a5}
\nccurve[angleA=180,angleB=240]{b1}{b2}
\endpspicture +
\pspicture[.4](-1,-1)(1,1)\hexanode
\ncline{a1}{b1}\ncline{a2}{b2}
\ncline{a5}{b5}\ncline{a6}{b6}
\ncarc[arcangle=6,doubleline=true]{a1}{a6}
\ncarc[arcangle=6]{a2}{a1}
\ncarc[arcangle=6]{a6}{a5}
\nccurve[angleA=300,angleB=0]{b3}{b4}
\endpspicture +
\pspicture[.4](-1,-1)(1,1)\hexanode
\ncline{a1}{b1}\ncline{a2}{b2}\ncline{a3}{b3}
\ncline{a4}{b4}
\ncarc[arcangle=6,doubleline=true]{a3}{a2}
\ncarc[arcangle=6]{a2}{a1}
\ncarc[arcangle=6]{a4}{a3}
\nccurve[angleA=60,angleB=120]{b5}{b6}
\endpspicture 
\right\}
+ \pspicture[.4](-1,-1)(1,1)\hexanode
\nccurve[angleA=180,angleB=240]{b1}{b2}
\nccurve[angleA=300,angleB=0]{b3}{b4}
\nccurve[angleA=60,angleB=120]{b5}{b6}
\endpspicture \label{hx}
\end{multline}
\caption{Hexagon relation}\label{relh}
\end{figure}

This presentation is not confluent. The reason is that the hexagon
relation is not a rewrite rule. If we have a hexagon in a diagram
then we can obtain a second diagram by interchanging the single
and double edges in the diagram. In order to make the hexagon relation
into a rewrite rule we need to know which of these two diagrams
is simpler.
\section{Centraliser algebras}\label{fpres}
In this section we define a sequence of finitely presented algebras.
For $n\ge 0$ the algebra $A_P(n)$ is generated by the set
\[ \{U_i,K_i,H_i\colon 1\le i\le n-1\} \]

The motivation for the relations is that they are satisfied in
$\mathcal{D}$. The diagrams for the generators are obtained from
\eqref{hom} by adding $i-1$ vertical lines on the left and
$n-i-1$ vertical lines on the right.
\begin{equation}\label{hom}
\psset{xunit=0.5,yunit=0.5}
\begin{array}{ccc}
U & K & H \\
\quad
\begin{pspicture}(2,-1)(4,1)
\psarc(3,1){0.3}{180}{360}
\psarc(3,-1){0.3}{0}{180}
\end{pspicture}\quad &
\quad
\begin{pspicture}(2,-1)(4,1)
\psline(2,1)(3,0.5)
\psline(2,-1)(3,-0.5)
\psline[doubleline=true](3,0.5)(3,-0.5)
\psline(3,0.5)(4,1)
\psline(3,-0.5)(4,-1)
\end{pspicture}\quad &
\quad
\begin{pspicture}(2,-1)(4,1)
\psline(2,1)(2.5,0)
\psline(2,-1)(2.5,0)
\psline[doubleline=true](2.5,0)(3.5,0)
\psline(3.5,0)(4,1)
\psline(3.5,0)(4,-1)
\end{pspicture}\quad
\end{array}
\psset{xunit=2,yunit=2}
\end{equation}

First we have the commuting relations. For $1\le i,j\le n-1$
with $|i-j|>1$, if $a\in\{U_i,K_i,H_i\}$ and  $b\in\{U_j,K_j,H_j\}$ 
then $ab=ba$.

The two string relations are:
\begin{equation*}
\begin{array}{ccc}
U_iU_i=\frac{[10][6][2]}{[5][3]} U_i & U_iH_i=[7]U_i & U_iK_i=0 \\
H_iU_i=[7]U_i & H_i^2 = [3]+\frac{[4]}{[2]}H_i+[4]U_i+[4]K_i & H_iK_i=-[5]K_i \\
K_iU_i=0 & K_iH_i=-[5]K_i & K_iK_i=\frac{[2][2][6]}{[3]}K_i
\end{array}
\end{equation*}
The eigenvalues of $H_i$ are $\{-1,[3],-[5],[7]\}$.
The relation for $H_i^2$ comes from the square relation \eqref{sq1}.

These relations imply that $A(2)$ is a commutative algebra with
basis $\{1,U_1,K_1,H_1\}$.

The three string relations are the following.
The first relations come from isotopy of diagrams.
\begin{eqnarray*}
U_iU_{i\pm 1}U_i &=& U_i \\
U_iH_{i\pm 1}U_i &=& 0 \\
U_iU_{i\pm 1}H_i &=& U_iK_{i\pm 1} \\
U_iU_{i\pm 1}K_i &=& U_iH_{i\pm 1} \\
H_iU_{i\pm 1}U_i &=& K_{i\pm 1}U_i \\
K_iU_{i\pm 1}U_i &=& H_{i\pm 1}U_i \\
K_iU_{i\pm 1}K_i &=& H_{i\pm 1}U_iH_{i\pm 1} 
\end{eqnarray*}
The next set of relations are relations in the ideal generated by
$U_i$ which come from the basic relations.
\begin{eqnarray*}
U_iK_{i\pm 1}U_i &=& [7]U_i \\
U_iK_{i\pm 1}K_i &=& -[5] U_iH_{i\pm 1} \\
K_iK_{i\pm 1}U_i &=& -[5] H_{i\pm 1}U_i \\
U_iH_{i\pm 1}U_i &=& 0 \\
U_iH_{i\pm 1}K_i &=& \frac{[2]^2[6]}{[3]} U_iH_{i\pm 1} \\
K_iH_{i\pm 1}U_i &=& \frac{[2]^2[6]}{[3]} H_{i\pm 1}U_i \\
U_iH_{i\pm 1}H_i &=& -[5] U_iH_{i\pm 1} \\
H_iH_{i\pm 1}U_i &=& -[5] H_{i\pm 1}U_i 
\end{eqnarray*}
The next relations are relations in the ideal generated by
$U_i$ and come from the square relations.
\begin{eqnarray*}
U_iK_{i\pm 1}H_i &=& [3] U_i + [4] U_iU_{i\pm 1} + \frac{[4]}{[2]} U_iK_{i\pm 1}
+ [4] U_iH_{i\pm 1} \\
H_iK_{i\pm 1}U_i &=& [3] U_i + [4] U_{i\pm 1}U_i + \frac{[4]}{[2]} K_{i\pm 1}U_i
+ [4] H_{i\pm 1}U_i 
\end{eqnarray*}

The next relation comes from one of the basic relations and the
other three come from the square relation \eqref{sq2}.
\begin{eqnarray*}
K_iH_{i\pm 1}K_i &=& 0 \\
K_iK_{i\pm 1}K_i &=& [2]^2K_i + [3]K_iU_{i\pm 1}K_i \\
H_iH_{i\pm 1}K_i &=& [2]^2U_{i\pm 1}K_i +[3]H_{i\pm 1}K_i \\
K_iH_{i\pm 1}H_i &=& [2]^2K_iU_{i\pm 1} +[3]K_iH_{i\pm 1} 
\end{eqnarray*}
The next relations come from the pentagon relation \eqref{pt}.
\begin{multline*}
K_iK_{i\pm 1}H_i=-\left( K_iK_{i\pm 1}+K_iU_{i\pm 1}+K_iU_{i\pm 1}H_i\right) \\
-[2]\left( K_i+K_iH_{i\pm 1}+K_iU_{i\pm 1}K_i+K_iU_{i\pm 1}U_i\right)
\end{multline*}
\begin{multline*}
H_iK_{i\pm 1}K_i=-\left( K_{i\pm 1}K_i+U_{i\pm 1}K_i+H_iU_{i\pm 1}K_i\right) \\
-[2]\left( K_i+H_{i\pm 1}K_i+K_iU_{i\pm 1}K_i+U_iU_{i\pm 1}K_i\right)
\end{multline*}
\begin{multline*}
H_iH_{i\pm 1}H_i=-\left(H_iH_{i\pm 1}+H_{i\pm 1}+H_{i\pm 1}H_i\right) \\
-[2]\left(U_iH_{i\pm 1}+H_{i\pm 1}U_i+K_iH_{i\pm 1}+H_{i\pm 1}K_i\right)
\end{multline*}
Finally we have the relation which comes from the hexagon relation \eqref{hx}.
\begin{multline*}
H_{i+1}\left( 1-[4]U_i - \frac{[4]}{[2]}H_i-[4]K_i \right) H_{i+1} = \\
H_i\left( 1-[4]U_{i+1} - \frac{[4]}{[2]}H_{i+1}-[4]K_{i+1} \right) H_i 
\end{multline*}
This can be rewritten using previous relations as
\begin{multline*}
H_{i+1}K_iH_{i+1}+H_{i+1}U_iH_{i+1}-U_{i+1}H_i-H_iU_{i+1}-K_{i+1}H_i-H_iK_{i+1}-U_{i+1}-K_{i+1} \\
= H_iK_{i+1}H_i+H_iU_{i+1}H_i-U_iH_{i+1}-H_{i+1}U_i-K_iH_{i+1}-H_{i+1}K_i-U_i-K_i 
\end{multline*}

The braid matrices are given by
\begin{eqnarray}\label{br} (q+q^{-1})\sigma_i &=&
-q^{-7}-q^{-2}U_i+q^{-4}K_i+q^{-5}H_i \\
(q+q^{-1})\sigma_i^{-1} &=&
-q^{7}-q^{2}U_i+q^{4}K_i+q^{5}H_i
\end{eqnarray}

\subsection{Yang-Baxter}
These relations can be described very succintly using the
Yang-Baxter equation (with spectral parameter). There are
two solutions of the Yang-Baxter equation. 

The first solution is discussed in \cite{MR1048745},\cite{MR1086739}
and \cite{MR1121816} and is associated with the affine quantum group
$U_q(B_3^{(1)})$. This solution can also be found using the
tensor product graph method introduced in \cite{MR92j:17013} and
\cite{MR92m:17031}. The tensor product graph is
\[ [0,0,0] \sr{-10} [0,1,0] \sr{-2} [0,0,2] \sr{6} [1,0,0] \]
The eigenvalues of $-(q-q^{-1})^4R_i(u)$ are
\begin{eqnarray*}
&[0,0,0]& (uq^5-u^{-1}q^{-5})(uq-u^{-1}q^{-1})(uq^{-3}-u^{-1}q^3) \\
&[0,1,0]& (u^{-1}q^5-uq^{-5})(uq-u^{-1}q^{-1})(uq^{-3}-u^{-1}q^3) \\
&[0,0,2]& (u^{-1}q^5-uq^{-5})(u^{-1}q-uq^{-1})(uq^{-3}-u^{-1}q^3) \\
&[1,0,0]& (u^{-1}q^5-uq^{-5})(u^{-1}q-uq^{-1})(u^{-1}q^{-3}-uq^3) 
\end{eqnarray*}
Then this satisfies
\begin{eqnarray*}
R_i(1) &=& [5][3] \\
R_i(q^5) &=& -[5][3]U_i \\
R_i(q^3) &=& -[3]K_i \\
R_i(q^2) &=& [3]H_i
\end{eqnarray*}
Hence $R_i(u)$ can be written as
\begin{multline*} (q-q^{-1})^2(q^2-q^{-2})R_i(u) = \\
-(uq^{-5}-u^{-1}q^5)(uq^{-3}-u^{-1}q^3)(uq^{-2}-u^{-1}q^2) \\
-(u-u^{-1})(uq^{-3}-u^{-1}q^3)(uq^{-2}-u^{-1}q^2) U_i \\
+(u-u^{-1})(uq^{-5}-u^{-1}q^5)(uq^{-2}-u^{-1}q^2) K_i \\
+(u-u^{-1})(uq^{-5}-u^{-1}q^5)(uq^{-3}-u^{-1}q^3) H_i \\
\end{multline*}
Then this satisfies the following relations.
\begin{eqnarray*}
R_i(u)R_j(v) &=& R_j(v)R_i(u) \text{  if $|i-j|>1$}\\
R_i(u)R_{i+1}(uv)R_i(v) &=& R_{i+1}(v)R_i(uv)R_{i+1}(u) \\
U_iU_{i\pm 1}R_i(u)&=& -U_iR_{i\pm 1}(q^5u^{-1}) \\
R_i(u)U_{i\pm 1}U_i &=& -R_{i\pm 1}(q^5u^{-1})U_i \\
R_i(u)U_{i+1}R_i(v) &=& -R_{i+1}(q^5u^{-1})U_iR_{i+1}(q^5v^{-1}) 
\end{eqnarray*}
\[ U_iR_{i\pm 1}(u)U_i =
-\frac{(uq^{-10}-u^{-1}q^{10})(uq^{-6}-u^{-1}q^{6})
(uq^{-2}-u^{-1}q^{2})}{(q-q^{-1})^3}  U_i \]

There is a second solution of the Yang-Baxter equation.
This solution is not explained by any known quantum group.
The solution is described by the tensor product graph
\[ [0,0,0] \sr{-6} [1,0,0] \sr{-4} [0,1,0] \sr{-2} [0,0,2] \]
The eigenvalues of $(q-q^{-1})^4S_i(u)$ are
\begin{eqnarray*}
&[0,0,0]& -(uq^3-u^{-1}q^{-3})(uq^2-u^{-1}q^{-2})(uq^{-1}-u^{-1}q) \\
&[1,0,0]& -(u^{-1}q^3-uq^{-3})(uq^2-u^{-1}q^{-2})(uq^{-1}-u^{-1}q) \\
&[0,1,0]& (u^{-1}q^3-uq^{-3})(u^{-1}q^2-uq^{-2})(uq^{-1}-u^{-1}q) \\
&[0,0,2]& (u^{-1}q^3-uq^{-3})(u^{-1}q^2-uq^{-2})(u^{-1}q-uq^{-1}) 
\end{eqnarray*}
These satisfy relations
\begin{eqnarray*}
S_i(u)S_j(v) &=& S_j(v)S_i(u) \text{  if $|i-j|>1$}\\
S_i(u)S_{i+1}(uv)S_i(v) &=& S_{i+1}(v)S_i(uv)S_{i+1}(u) \\
\end{eqnarray*}

The sequence of algebras $A_P(n)$ can also be defined by the commuting
relations, the two string relations and both Yang-Baxter equations.
\subsection{Representations}\label{rep}
Next we give irreducible representations of dimensions 1,2,3 and
4 of $A(3)$. In each case we only give the
matrices representing $U_1$, $K_1$,$H_1$ and $\sigma_1$. The
reason is that, for each of these representations of dimension $n$,
the matrices representing $U_2$, $K_2$, $H_2$ and $\sigma_2$
are obtained by applying the involution
\[ A_{ij} \leftrightarrow A_{n-j+1,n-i+1} \]
Also the matrices for $\sigma_i^{-1}$ are given by applying the
involution $q\leftrightarrow q^{-1}$ to the entries of the matrices
for $\sigma_1$.

The reason we have included the matrices for $\sigma_1$ even
though they are determined by \eqref{br} is that we found these
matrices first using the results in \cite{MR1815266} and then
we calculated the matrices representing the generators from
these.

The four dimensional representation is given by
\[
U_1 =
\left[ \begin {array}{cccc} \frac{[10][6][2]}{[5][3]}
 &-\frac{[10][6]}{[5][3]}
 &-\frac{[6][4]}{[3][2]}
  &1\\\noalign{\medskip}0&0&0&0\\\noalign{\medskip}0&0&0&0
  \\\noalign{\medskip}0&0&0&0\end {array} \right] 
\]
%\[ U_2 =
%\left[ \begin {array}{cccc} 0&0&0&0\\\noalign{\medskip}0&0&0&0
%\\\noalign{\medskip}0&0&0&0\\\noalign{\medskip}1&
%-\frac{[6][4]}{3][2]}&
%-\frac{[10][6]}{[5][3]}&
%\frac{[10][6][2]}{[5][3]}
   %\end {array} \right] \]
\[ K_1 =
\left[ \begin {array}{cccc}
0&[2]\frac{[6]}{[3]} &-\frac {[6]}{[3]}&-1
\\\noalign{\medskip}0&[2]^2\frac{[6]}{[3]}&-[2]\frac{[6]}{[3]}
&-[2]\\\noalign{\medskip}0&0&0&0\\\noalign{\medskip}0&0&0&0\end {array}
\right] 
\]
%\[ K_2 =
%\left[ \begin {array}{cccc} 0&0&0&0\\\noalign{\medskip}0&0&0&0
%\\\noalign{\medskip}-{\frac {{q}^{2}+1}{q}}&-{\frac { \left( {q}^{2}+1
 %\right)  \left( {q}^{6}+1 \right) }{{q}^{4}}}&{\frac { \left( {q}^{2}
 %+1 \right) ^{2} \left( {q}^{6}+1 \right) }{{q}^{5}}}&0
 %\\\noalign{\medskip}-1&-{\frac { \left( {q}^{6}+1 \right)  
 %}{{q}^{3}}}&{\frac { \left( 1+{q}^{6}
  %\right)  \left( {q}^{2}+1 \right) }{{q}^{4}}}&0\end {array}
   %\right] 
%\]
\[ H_1 =
\left[ \begin {array}{cccc} [7]&-[6]&0&0\\\noalign{\medskip}0
 &-[5]&[4]&0\\\noalign{\medskip}0&0&[3]&-[2]\\\noalign{\medskip}0&0&0
 &-1\end {array} \right] 
\]
%\[ H_2 =
%\left[ \begin {array}{cccc} -1&0&0&0\\\noalign{\medskip}-{\frac {{q}^
%{2}+1}{q}}&{\frac {{q}^{4}+{q}^{2}+1}{{q}^{2}}}&0&0
%\\\noalign{\medskip}0&{\frac {{q}^{6}+{q}^{4}+{q}^{2}+1}{{q}^{3}}}&-{
%\frac {{q}^{8}+{q}^{6}+{q}^{4}+{q}^{2}+1}{{q}^{4}}}&0
%\\\noalign{\medskip}0&0&-{\frac { \left( {q}^{6}+1 \right)  
 %\left( {q}^{4}+{q}^{2}+1 \right) }{{q}^{5}}}&{
%\frac {{q}^{6}+{q}^{8}+{q}^{4}+{q}^{10}+{q}^{12}+{q}^{2}+1}{{q}^{6}}}
%\end {array} \right]
%\]
\[ \sigma_1 =
\left[ \begin {array}{cccc} -{q}^{6}&q-{q}^{3}+{q}^{5}&{q}^{-2}-1+{q}
^{2}&-{q}^{-3}\\\noalign{\medskip}0&1&{q}^{-3}-{q}^{-1}&-{q}^{-4}
\\\noalign{\medskip}0&0&{q}^{-4}&-{q}^{-5}\\\noalign{\medskip}0&0&0&-{
q}^{-6}\end {array} \right]
\]

For $1\le k\le 3$ there is a representation of dimension $4-k$.
The matrices representing $U_1$, $K_1$, $H_1$ and $\sigma_1^{\pm 1}$
are obtained by deleting the first $k$ rows and columns from the
matrix representing the same element in the four dimensional
representation.

This means that there is a spanning set 
for $A_P(3)$ which consists of the words of length
at most two the hexagon $H_1K_2H_1$ and the following four words
in the ideal generated by $U_2$.
\[ \begin{array}{cc}
H_1U_2H_1 & H_1U_2K_1 \\ K_1U_2H_1 & K_1U_2K_1
\end{array} \]
This gives a total of thirty words in the spanning set.
The sum of the squares of the dimensions of the above
representations is thirty and these thirty words are linearly
independent in this thirty dimensional algebra. This shows
that these thirty words are a basis for $A(3)$.

It follows from these relations that the criterion of
\cite{MR1465031} is satisfied. This means that 
\begin{multline}\label{w}
A_P(n+1) = \\
A_P(n) + A_P(n)U_nA_P(n) + A_P(n)K_nA_P(n) + A_P(n)H_nA_P(n)
\end{multline}
The first consequence of this is that there is a conditional
expectation $\varepsilon_n\colon A_P(n)\rightarrow A_P(n-1)$.
This satisfies $U_naU_n = \varepsilon(a)U_n$ for all
$a\in A_P(n)$ and is uniquely determined by the conditions
\begin{eqnarray*}
\varepsilon(a) &=& \frac{[10][6][2]}{[5][3]} a \\
\varepsilon(aU_{n-1}a^\prime) &=&  aa^\prime \\
\varepsilon(aK_{n-1}a^\prime) &=&  aa^\prime \\
\varepsilon(aH_{n-1}a^\prime) &=&  [7]aa^\prime
\end{eqnarray*}
for all $a,a^\prime\in A_P(n-2)$. These conditional expectations
then determine a trace map $\tau_n$ on $A_P(n)$ for each $n>0$
by $\tau_{n+1}(a)=\frac{[5][3]}{[10][6][2]}
\tau_n(\varepsilon_{n+1}(a))$ for all
$a\in A_P(n)$.

The second consequence of \eqref{w} is that for any $n\ge 0$,
there is a finite set of words in the generators which span $A_P(n)$.
In fact we can be much more specific.

\begin{defn}\label{irr}
Take an array of integers $\{ a(i,j)| i\ge 1, j\ge 1, i+j\le n+1\}$
such that $0\le a(i,j)\le 3$,
$a(i,j)\le a(i^\prime,j)$ if $i<i^\prime$,
$a(i,j)\le a(i,j^\prime)$ if $j<j^\prime$.
Then associated to this array is the following word in 
non-commuting indeterminates $E_i^{(k)}$ where $1\le i\le n-1$
and $0\le k\le 3$
\[ \prod_{j=n-1}^1\left(\prod_{i=1}^{n-j}E_{n-i-j+1}^{a(i,j)}
\right) \]
Then substitute $1$ for $E_i^{(0)}$, $U_i$ for $E_i^{(1)}$,
$K_i$ for $E_i^{(2)}$, $K_i$ for $E_i^{(3)}$ to get a word
in the generators of $A_P(n)$. We will call these the
irreducible words.
\end{defn}

It is clear that the irreducible words span $A_P(n)$ since if we take
any irreducible word and multiply by a generator then by
using the relations we can write this word as a linear 
combination of irreducible words. In fact the irreducible words are
a basis and we will give two proofs of this; one in \S \ref{conc}
using representation theory and one in \S \ref{conf} using the
diamond lemma.

The third consequence of \eqref{w} is:
\begin{lem}\label{l}
For all $n>0$, $U_nA_P(n+1)U_n=A_P(n-1)U_n$.
\end{lem}
\begin{pf} Let $a\in A_P(n+1)$. Then by \eqref{w}, $a$
can be written as a linear combination of terms of the form
$bX_nb^\prime$ where $b,b^\prime \in A_P(n)$
and $X_n\in \{1,U_n,K_n,H_n\}$. Then applying \eqref{w} a second
time, each $b$ can be written as a linear combination of elements
of the form $cY_{n-1}c^\prime$ where
$c,c^\prime \in A_P(n-1)$ and $Y_{n-1}\in \{1,U_{n-1},K_{n-1},H_{n-1}\}$.
Hence $U_naU_n$ can be written as a linear combination of terms
of the form
\[ cU_nY_{n-1}X_n(c^\prime b^\prime)U_n \]
where $c\in A_P(n-1)$ and $(c^\prime b^\prime)\in A_P(n)$.
Each word of the form $U_nY_{n-1}X_n$ can be written as a linear
combination of the words 
\[ \{U_n,U_nU_{n-1},U_nK_{n-1},U_nH_{n-1}\} \]
using the defining relations.
This shows that $U_naU_n$ is a element of $U_nA_P(n)U_n$ and so
is an element of $A_P(n-1)U_n$.
\end{pf}
\section{Crystal graphs}\label{crystal}
Next we discuss the combinatorics of Littelman paths for the
spin representation. The first observation is that the weights of this
representation are the orbit of the highest weight under the action
of the Weyl group. This implies that all weight spaces have dimension
one and that the Littelman paths are just
the straight lines from the origin to the weight. Then when we 
concatenate these paths we get a sequence of weights 
$\omega_1,\omega_2,\ldots $ such that for all $i>2$,
$\omega_{i+1}-\omega_i$ is a weight of the representation.
Furthermore the sequence of weights corresponds to a dominant path
if and only if every weight in the sequence is dominant.

Next instead of working with the basis consisting of the
fundamental weights we change basis by
\begin{eqnarray*}
(w_1,w_2,w_3) &\mapsto& (2w_1+2w_2+w_3,2w_2+w_3,w_3) \\
(s_1,s_2,s_3) &\mapsto& (\frac{s_1-s_2}{2},\frac{s_2-s_3}{2},s_3)
\end{eqnarray*}
In the new coordinates the weights of the representation are
the eight vectors $(\pm 1,\pm 1,\pm 1)$. Also the dominant weights
correspond to vectors $(s_1,s_2,s_3)$ such that 
$s_1\ge s_2\ge s_3\ge 0$. This shows that Littelman paths correspond
to triples of non-crossing Dyck paths.

The number of triples of non-crossing Dyck paths is calculated in
\cite{MR927758} where it is shown to be
\[ \prod_{1\le i\le j\le n}\frac{i+j+6}{i+j} \]
This sequence appears in the On-Line Encyclopedia of Integer
Sequences as sequence A006149.

%Let $\Delta$ be the spin representation. There are no non-zero
%invariant tensors in $\otimes^n\Delta$ for $n$ odd. Let $a(k)$
%be the dimension of the vector space of invariant tensors
%in $\otimes^{2k}\Delta$ for $k\ge 0$. This sequence starts
%\begin{center}
%\begin{tabular}{rrrrrrrr}
%0 & 1 & 2 & 3 & 4 & 5 & 6 & 7 \\
%\hline
%1 & 1 & 4 & 30 & 330 & 4719 & 81796 & 1643356 \\
%\end{tabular}
%\end{center}
%This sequence appears in the On-Line Encyclopedia of Integer
%Sequences as sequence A006149 and is studied in \cite{mstc}.
%The sequence is determined by
%the initial condition $a(0)=1$ and the recurrence relation
%\[
%a(k+1)=\frac{(2k+1)(2k+3)(2k+5)}{(k+4)(k+5)(k+6)}8a(k)
%=\frac{(2k+6)!}{(2k)!}\frac{k!}{(k+6)!}a(k)
%\]
%Equivalently, the ordinary generating function is the
%hypergeometric function
%\[ {}_4F_3([1,1/2,3/2,5/2];[4,5,6];64x) \]
%
%This gives the following expression for $a(k)$
%\[ a(k) = \frac%
%{3!5!(2k+4)!(2k+2)!(2k)!}%
%{(k+6)!(k+5)!(k+4)!(k+3)!(k+2)!(k+1)!} \]

%This generalises to $B_n$ as follows
%\[ a(k+1)=\frac{\prod_{i=1}^n(2k+2i-1)}%
%{\prod_{i=1}^n (k+n+i)}2^n a(k) \]
%The case $n=4$ is sequence A006150 and the case $n=5$ is
%sequence A006151. The case $n=1$ gives the Catalan numbers.
%Equivalently, the ordinary generating function is the
%hypergeometric function
%\[ {}_{n+1}F_n([1,1/2,3/2,\ldots ,(2n-1)/2];[n+1,n+2,\ldots ,2n];2^{2n}x) \]

In this section we use the algorithms introduced in
\cite{math.CO/0507112} to find a set of diagrams which gives a basis
for the invariant tensors.

\begin{defn}\label{cut} Assume we are given a diagram.
Let $A$ and $B$ be two boundary points which are not marked points.
Then a cut path from $A$ to $B$ is a path from $A$ to $B$ such that
each component of the intersection with the embedded graph is either
an isolated transverse intersection point or else is an edge of the
graph.
\end{defn}
The diagrams for these four cases and the associated weights are:
\begin{center}
\begin{pspicture}(-2,-1.5)(7,0.5)
\psline[linestyle=dashed](-2,0)(-1,0)
\psline(-1.5,0.5)(-1.5,-0.5)
\psline[linestyle=dashed](0,0)(1,0)
\psline[doubleline=true](0.5,0.5)(0.5,-0.5)
\psline[linestyle=dashed](2,0)(4,0)
\psline(2,0.5)(2.5,0)
\psline(2,-0.5)(2.5,0)
\psline[doubleline=true](2.5,0)(3.5,0)
\psline(3.5,0)(4,0.5)
\psline(3.5,0)(4,-0.5)
\psline[linestyle=dashed](5,0)(7,0)
\psline[doubleline=true](5,0.5)(5.5,0)
\psline(5,-0.5)(5.5,0)
\psline(5.5,0)(6.5,0)
\psline[doubleline=true](6.5,0)(7,0.5)
\psline(6.5,0)(7,-0.5)
\rput(-1.5,-1){$(0,0,1)$}
\rput(0.5,-1){$(1,0,0)$}
\rput(3,-1){$(0,1,0)$}
\rput(6,-1){$(0,1,0)$}
\end{pspicture}
\end{center}
The weight of a cut path is $(a,b,c)$ if it crosses $a$ double
edges transversally, $c$ single edges transversally and contains
$b$ edges. A cut path is minimal if there is no cut path with the
same endpoints and lower weight.

Now suppose we have a diagram with $n$ boundary points which are
all single edges. Then we draw the diagram in a triangle $AXY$
with the edge $XY$ horizontal and with all boundary points on the
edge $XY$. Then we choose marked points $X_0=X,X_1,\ldots X_n=Y$
such there is precisely one boundary point between each marked
point. Now for $0\le i\le n$ let $\omega_i$ be the minimum weight
such there is a cut path from $A$ to $X_i$ of weight $\omega_i$.
Then $\omega_0=(0,0,0)=\omega_n$.

\begin{figure}
\psset{xunit=0.07cm,yunit=0.07cm}
\psset{yunit=1.732}
\newcommand{\TVA}{
\psset{linewidth=1pt,linestyle=dotted}
\pcline(-12,4)(12,4)
\pcline(12,4)(0,-8)\Aput{(1,0,0)}
\pcline(-12,4)(0,-8)\Bput{(0,0,0)}
\psset{linewidth=2pt,linestyle=solid}
\psline[doubleline=true,linewidth=1pt](0,4)(0,0)
\psline[doubleline=true,linewidth=1pt](0,0)(6,-2) }

\newcommand{\TVB}{
\psset{linewidth=1pt,linestyle=dotted}
\pcline(-12,4)(12,4)
\pcline(12,4)(0,-8)\Aput{(0,1,0)}
\pcline(-12,4)(0,-8)\Bput{(1,0,0)}
\psset{linewidth=2pt,linestyle=solid}
\psline[doubleline=true,linewidth=1pt](6,4)(6,2)
\psline(0,0)(6,2)
\psline(0,0)(0,-4)
\psline[doubleline=true,linewidth=1pt](-3,-5)(0,-4)
\psline(6,2)(9,1)
\psline(0,-4)(3,-5) }

\newcommand{\TVC}{
\psset{linewidth=1pt,linestyle=dotted}
\pcline(-12,4)(12,4)
\pcline(12,4)(0,-8)\Aput{(0,0,2)}
\pcline(-12,4)(0,-8)\Bput{(0,1,0)}
\psset{linewidth=2pt,linestyle=solid}
\psline[doubleline=true,linewidth=1pt](6,4)(6,2)
\psline(-9,1)(-6,2)
\psline(-6,2)(0,0)
\psline(0,0)(6,2)
\psline(6,2)(9,1)
\psline[doubleline=true,linewidth=1pt](0,0)(0,-4)
\psline(0,-4)(3,-5)
\psline(0,-4)(-3,-5) }

\newcommand{\TVD}{
\psset{linewidth=1pt,linestyle=dotted}
\pcline(-12,4)(12,4)
\pcline(12,4)(0,-8)\Aput{(0,0,1)}
\pcline(-12,4)(0,-8)\Bput{(0,0,1)}
\psset{linewidth=2pt,linestyle=solid}
\psline[doubleline=true,linewidth=1pt](0,4)(0,0)
\psline(0,0)(-6,-2) 
\psline(0,0)(6,-2) }

\newcommand{\TVE}{
\psset{linewidth=1pt,linestyle=dotted}
\pcline(-12,4)(12,4)
\pcline(12,4)(0,-8)\Aput{(0,1,0)}
\pcline(-12,4)(0,-8)\Bput{(0,0,2)}
\psset{linewidth=2pt,linestyle=solid}
\psline[doubleline=true,linewidth=1pt](-6,4)(-6,2)
\psline(-9,1)(-6,2)
\psline(-6,2)(0,0)
\psline(0,0)(6,2)
\psline(6,2)(9,1)
\psline[doubleline=true,linewidth=1pt](0,0)(0,-4)
\psline(0,-4)(3,-5)
\psline(0,-4)(-3,-5) }

\newcommand{\TVF}{
\psset{linewidth=1pt,linestyle=dotted}
\pcline(-12,4)(12,4)
\pcline(12,4)(0,-8)\Aput{(1,0,0)}
\pcline(-12,4)(0,-8)\Bput{(0,1,0)}
\psset{linewidth=2pt,linestyle=solid}
\psline[doubleline=true,linewidth=1pt](-6,4)(-6,2)
\psline(0,0)(-6,2)
\psline(0,0)(0,-4)
\psline[doubleline=true,linewidth=1pt](3,-5)(0,-4)
\psline(-6,2)(-9,1)
\psline(0,-4)(-3,-5) }

\newcommand{\TVG}{
\psset{linewidth=1pt,linestyle=dotted}
\pcline(-12,4)(12,4)
\pcline(12,4)(0,-8)\Aput{(0,0,0)}
\pcline(-12,4)(0,-8)\Bput{(1,0,0)}
\psset{linewidth=2pt,linestyle=solid}
\psline[doubleline=true,linewidth=1pt](0,4)(0,0)
\psline[doubleline=true,linewidth=1pt](0,0)(-6,-2) }

\begin{center}
\begin{pspicture}(0,-10)(36,170)
\rput(0,180){\rnode{v1}{\TVA}}
\rput(0,150){\rnode{v2}{\TVB}}
\rput(0,120){\rnode{v3}{\TVC}}
\rput(0,90){\rnode{v4}{\TVD}}
\rput(0,60){\rnode{v5}{\TVE}}
\rput(0,30){\rnode{v6}{\TVF}}
\rput(0,0){\rnode{v7}{\TVG}}
\psset{nodesep=30pt,arrows=->,arrowsize=3pt 3,linewidth=1pt}
\ncline{v1}{v2}\mput*{1}
\ncline{v2}{v3}\mput*{2}
\ncline{v3}{v4}\mput*{3}
\ncline{v4}{v5}\mput*{3}
\ncline{v5}{v6}\mput*{2}
\ncline{v6}{v7}\mput*{1}
\end{pspicture}
\end{center}
\caption{Crystal graph for the vector representation}\label{vec}
\end{figure}
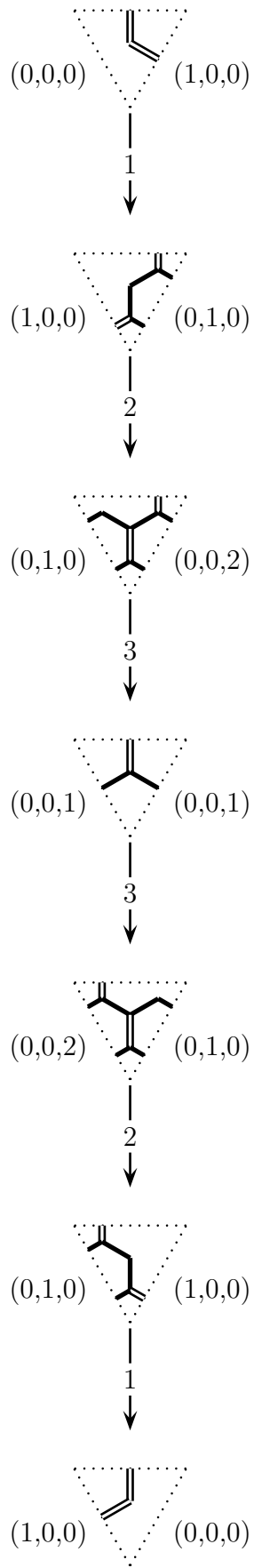

The diagram model for the crystal graph of the eight dimensional spin
representation is given in Figure \ref{sp} and
the diagram model for the crystal graph of the seven dimensional vector
representation is given in Figure \ref{vec}.
These crystal graphs are given in \cite[\S 8.1]{MR1881971}.
\begin{rem}
There is a different interpretation of these labelled directed graphs.
These two representations are miniscule which means
that the action of the Weyl group $W$ on the weights is transitive.
The stabiliser of a point is a parabolic subgroup $W_0$. Each coset
has a unique representative of minimal length. In these
examples $W$ has type $B_3$; for the vector representation
$W_0$ has type $B_2$ and for the spin representation $W_0$ has type
$A_2$. The Weyl group of type
$B_3$ is generated by the reflections in the three simple roots,
$s_1$,$s_2$,$s_3$. For each vertex, take a directed path from the
highest weight (at the top) to the vertex. Then take the sequence of
labels and regard it as a word in the generators. These words are
then reduced expressions for the minimal length coset representative.
\end{rem}
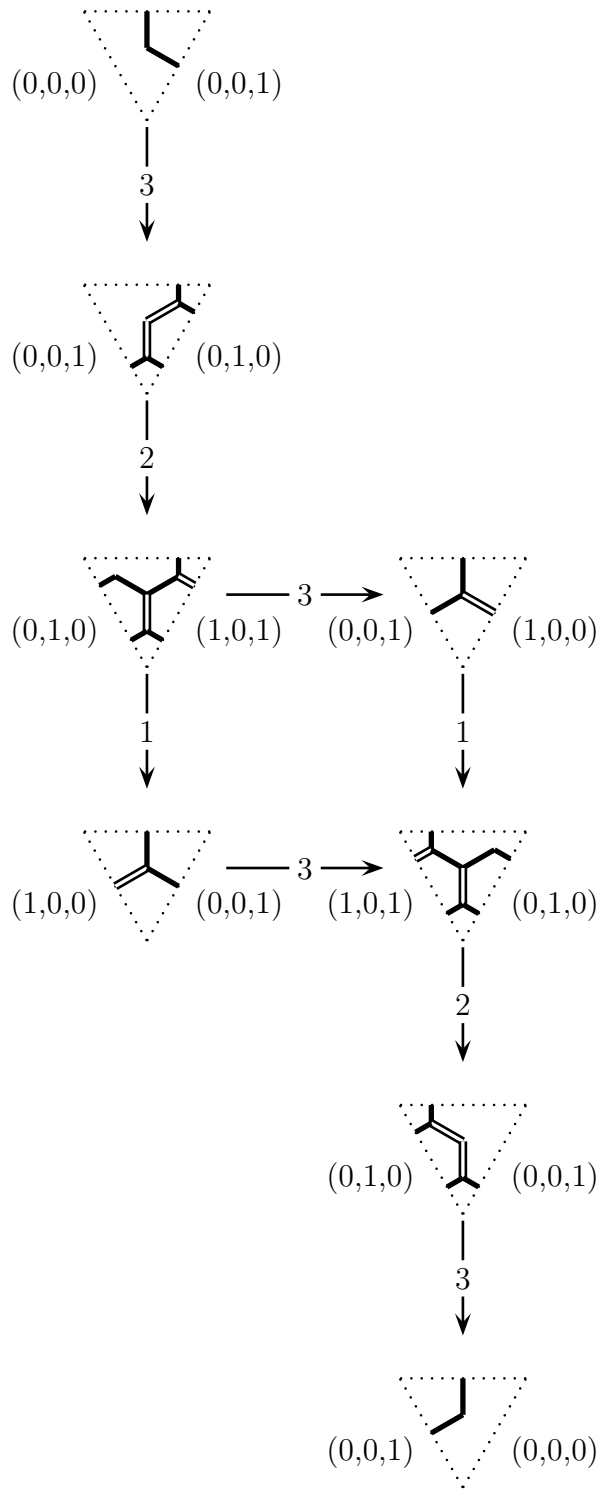
\begin{figure}
\psset{xunit=0.07cm,yunit=0.07cm}
\newcommand{\TBA}{
\psset{linewidth=1pt,linestyle=dotted}
\pcline(-12,4)(12,4)
\pcline(12,4)(0,-8)\Aput{(0,0,1)}
\pcline(-12,4)(0,-8)\Bput{(0,0,0)}
\psset{linewidth=2pt,linestyle=solid}
\psline(0,4)(0,0)
\psline(0,0)(6,-2) }

\newcommand{\TBD}{
\psset{linewidth=1pt,linestyle=dotted}
\pcline(-12,4)(12,4)
\pcline(12,4)(0,-8)\Aput{(0,0,1)}
\pcline(-12,4)(0,-8)\Bput{(1,0,0)}
\psset{linewidth=2pt,linestyle=solid}
\psline(0,4)(0,0)
\psline[doubleline=true,linewidth=1pt](0,0)(-6,-2) 
\psline(0,0)(6,-2) }

\newcommand{\TBE}{
\psset{linewidth=1pt,linestyle=dotted}
\pcline(-12,4)(12,4)
\pcline(12,4)(0,-8)\Aput{(1,0,0)}
\pcline(-12,4)(0,-8)\Bput{(0,0,1)}
\psset{linewidth=2pt,linestyle=solid}
\psline(0,4)(0,0)
\psline(0,0)(-6,-2) 
\psline[doubleline=true,linewidth=1pt](0,0)(6,-2) }

\newcommand{\TBH}{
\psset{linewidth=1pt,linestyle=dotted}
\pcline(-12,4)(12,4)
\pcline(12,4)(0,-8)\Aput{(0,0,0)}
\pcline(-12,4)(0,-8)\Bput{(0,0,1)}
\psset{linewidth=2pt,linestyle=solid}
\psline(0,4)(0,0)
\psline(0,0)(-6,-2) }

\newcommand{\TGC}{
\psset{linewidth=1pt,linestyle=dotted}
\pcline(-12,4)(12,4)
\pcline(12,4)(0,-8)\Aput{(0,1,0)}
\pcline(-12,4)(0,-8)\Bput{(1,0,1)}
\psset{linewidth=2pt,linestyle=solid}
\psline(-6,4)(-6,2)
\psline[doubleline=true,linewidth=1pt](-9,1)(-6,2)
\psline(-6,2)(0,0)
\psline(0,0)(6,2)
\psline(6,2)(9,1)
\psline[doubleline=true,linewidth=1pt](0,0)(0,-4)
\psline(0,-4)(3,-5)
\psline(0,-4)(-3,-5) }

\newcommand{\TGE}{
\psset{linewidth=1pt,linestyle=dotted}
\pcline(-12,4)(12,4)
\pcline(12,4)(0,-8)\Aput{(1,0,1)}
\pcline(-12,4)(0,-8)\Bput{(0,1,0)}
\psset{linewidth=2pt,linestyle=solid}
\psline(6,4)(6,2)
\psline(-9,1)(-6,2)
\psline(-6,2)(0,0)
\psline(0,0)(6,2)
\psline[doubleline=true,linewidth=1pt](6,2)(9,1)
\psline[doubleline=true,linewidth=1pt](0,0)(0,-4)
\psline(0,-4)(3,-5)
\psline(0,-4)(-3,-5) }

\newcommand{\TGF}{
\psset{linewidth=1pt,linestyle=dotted}
\pcline(-12,4)(12,4)
\pcline(12,4)(0,-8)\Aput{(0,0,1)}
\pcline(-12,4)(0,-8)\Bput{(0,1,0)}
\psset{linewidth=2pt,linestyle=solid}
\psline(-6,4)(-6,2)
\psline[doubleline=true,linewidth=1pt](0,0)(-6,2)
\psline[doubleline=true,linewidth=1pt](0,0)(0,-4)
\psline(3,-5)(0,-4)
\psline(-6,2)(-9,1)
\psline(0,-4)(-3,-5) }

\newcommand{\TGB}{
\psset{linewidth=1pt,linestyle=dotted}
\pcline(-12,4)(12,4)
\pcline(12,4)(0,-8)\Aput{(0,1,0)}
\pcline(-12,4)(0,-8)\Bput{(0,0,1)}
\psset{linewidth=2pt,linestyle=solid}
\psline(6,4)(6,2)
\psline[doubleline=true,linewidth=1pt](0,0)(6,2)
\psline[doubleline=true,linewidth=1pt](0,0)(0,-4)
\psline(-3,-5)(0,-4)
\psline(6,2)(9,1)
\psline(0,-4)(3,-5) }

\begin{center}
\psset{yunit=1.732}
\begin{pspicture}(0,-10)(72,169)
\rput(0,150){\rnode{v1}{\TBA}}
\rput(0,120){\rnode{v2}{\TGB}}
\rput(0,90){\rnode{v3}{\TGE}}
\rput(0,60){\rnode{v4}{\TBD}}
\rput(60,90){\rnode{v5}{\TBE}}
\rput(60,60){\rnode{v6}{\TGC}}
\rput(60,30){\rnode{v7}{\TGF}}
\rput(60, 0){\rnode{v8}{\TBH}}
\psset{nodesep=30pt,arrows=->,arrowsize=3pt 3,linewidth=1pt}
\ncline{v1}{v2}\mput*{3}
\ncline{v2}{v3}\mput*{2}
\ncline{v3}{v4}\mput*{1}
\ncline{v3}{v5}\mput*{3}
\ncline{v4}{v6}\mput*{3}
\ncline{v5}{v6}\mput*{1}
\ncline{v6}{v7}\mput*{2}
\ncline{v7}{v8}\mput*{3}
\end{pspicture}
\end{center}
\caption{Crystal graph for the spin representation}\label{sp}
\end{figure}

Next we discuss the tensor product rule. This tensor product rule
gives a procedure which associates a triangular diagram to any
word in the weights of the spin representation and the vector
representation. This tensor product rule is similar to the tensor
product rule given in \cite{math.CO/0507112}. First we draw a 
triangular grid. The triangles in top row is the sequence of triangles
associated to the word. This also gives the weights of the edges 
of these triangles. Next we label each edge of the triangular grid
by a dominant weight using the rule below. This rule is derived from the
tensor product of crystals.
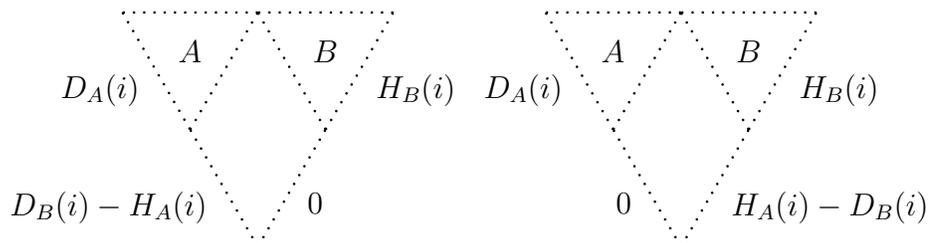
\psset{xunit=0.10cm,yunit=0.1732cm}
\psset{xunit=0.75,yunit=0.75}
\begin{center}
\begin{pspicture}(-12,-20)(111,4)

\rput(0,0){
\psset{linewidth=1pt,linestyle=dotted}
\pcline(-12,4)(12,4)
\pcline(12,4)(0,-8)
\pcline(0,-8)(-12,4)\Aput{$D_A(i)$}

\pcline(12,4)(36,4)
\pcline(36,4)(24,-8)\Aput{$H_B(i)$}
\pcline(12,4)(24,-8)

\pcline(0,-8)(12,-20)\Bput{$D_B(i)-H_A(i)$}
\pcline(24,-8)(12,-20)\Aput{$0$}
}
\rput(0,0){$A$}
\rput(24,0){$B$}

\rput(75,0){
\psset{linewidth=1pt,linestyle=dotted}
\pcline(-12,4)(12,4)
\pcline(12,4)(0,-8)
\pcline(-12,4)(0,-8)\Bput{$D_A(i)$}

\pcline(12,4)(36,4)
\pcline(36,4)(24,-8)\Aput{$H_B(i)$}
\pcline(12,4)(24,-8)

\pcline(0,-8)(12,-20)\Bput{$0$}
\pcline(24,-8)(12,-20)\Aput{$H_A(i)-D_B(i)$}
\rput(0,0){$A$}
\rput(24,0){$B$}
}
\end{pspicture}
\end{center}
\psset{xunit=1cm,yunit=1cm}

Now we complete the triangular diagram by drawing a graph
in each diamond. Note that the weights $H_A$
and $D_B$ are both elements of the set 
\[ \{(0,0,0),(0,0,1),(0,1,0),(1,0,0),(1,0,1),(0,0,2)\} \]
This gives thirty six different diamonds. Furthermore note
that for each of these diamonds the other two
weights are also elements of this set. This implies that
for any word every diamond will have one of these thirty six
labellings. Therefore to complete the diagram it is sufficient
to know how to fill in a diamond with each of these labellings.

The general principle for filling in these diagrams is that 
each edge of the diamond is required to be a minimal cut path
whose weight is given by the label and the diagram is required
to be irreducible. By construction the two paths from the bottom
of the diamond to the top are cut paths of the same weight.

The simplest case is where two of the edges are
labelled with the zero weight. These can be filled in
using the following diagrams:
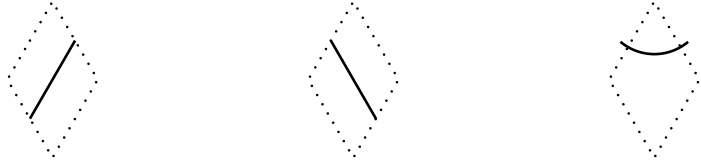
\begin{center}
\psset{xunit=0.1cm,yunit=0.1732cm}
\begin{pspicture}(-6,-8)(86,4)

\rput(0,0){
\psset{linewidth=1pt,linestyle=dotted}
\pspolygon(0,4)(6,-2)(0,-8)(-6,-2)
\psset{linestyle=solid}
\psline(3,1)(-3,-5)
}
\rput(40,0){
\psset{linewidth=1pt,linestyle=dotted}
\pspolygon(0,4)(6,-2)(0,-8)(-6,-2)
\psset{linestyle=solid}
\psline(-3,1)(3,-5)
}
\rput(80,0){
\psset{linewidth=1pt,linestyle=dotted}
\pspolygon(0,4)(6,-2)(0,-8)(-6,-2)
\psset{linestyle=solid}
\psarc(0,4){0.7}{-130}{-50}
}
\end{pspicture}
\psset{xunit=1cm,yunit=1cm}
\end{center}
In these diagrams any solid line can be removed or replaced
by any sequence of parallel lines.

This leaves twenty cases. However these come in pairs since
any diamond can be reflected in a vertical line. This leaves
ten cases which we describe below.

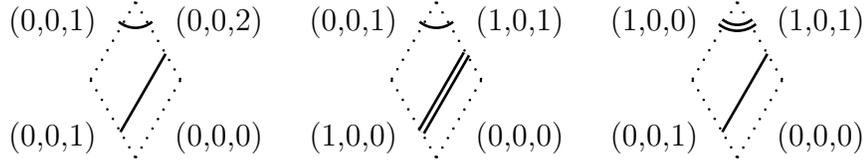
\begin{center}
\psset{xunit=0.1cm,yunit=0.1732cm}
\begin{pspicture}(-6,-8)(86,4)
\rput(0,0){
\psset{linewidth=1pt,linestyle=dotted}
\pcline(0,4)(6,-2)\Aput{(0,0,2)}
\pcline(6,-2)(0,-8)\Aput{(0,0,0)}
\pcline(0,-8)(-6,-2)\Aput{(0,0,1)}
\pcline(-6,-2)(0,4)\Aput{(0,0,1)}
\psset{linestyle=solid}
\psline(4,0)(-2,-6)
\psarc(0,4){0.3464}{-130}{-50}
}
\rput(40,0){
\psset{linewidth=1pt,linestyle=dotted}
\pcline(0,4)(6,-2)\Aput{(1,0,1)}
\pcline(6,-2)(0,-8)\Aput{(0,0,0)}
\pcline(0,-8)(-6,-2)\Aput{(1,0,0)}
\pcline(-6,-2)(0,4)\Aput{(0,0,1)}
\psset{linestyle=solid}
\psline[doubleline=true,linewidth=1pt](4,0)(-2,-6)
\psarc(0,4){0.3464}{-130}{-50}
}
\rput(80,0){
\psset{linewidth=1pt,linestyle=dotted}
\pcline(0,4)(6,-2)\Aput{(1,0,1)}
\pcline(6,-2)(0,-8)\Aput{(0,0,0)}
\pcline(0,-8)(-6,-2)\Aput{(0,0,1)}
\pcline(-6,-2)(0,4)\Aput{(1,0,0)}
\psset{linestyle=solid}
\psline(4,0)(-2,-6)
\psarc[doubleline=true,linewidth=1pt](0,4){0.3464}{-130}{-50}
}
\end{pspicture}
\psset{xunit=1cm,yunit=1cm}
\end{center}

Here are three more cases.
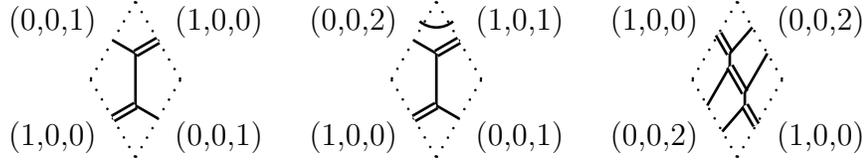
\begin{center}
\psset{xunit=0.1cm,yunit=0.1732cm}
\begin{pspicture}(-6,-8)(86,4)
\rput(0,0){
\psset{linewidth=1pt,linestyle=dotted}
\pcline(0,4)(6,-2)\Aput{(1,0,0)}
\pcline(6,-2)(0,-8)\Aput{(0,0,1)}
\pcline(0,-8)(-6,-2)\Aput{(1,0,0)}
\pcline(-6,-2)(0,4)\Aput{(0,0,1)}
\psset{linestyle=solid}
\psline(-3,1)(0,0)
\psline[doubleline=true,linewidth=1pt](3,1)(0,0)
\psline(0,0)(0,-4)
\psline[doubleline=true,linewidth=1pt](0,-4)(-3,-5)
\psline(0,-4)(3,-5)
}
\rput(40,0){
\psset{linewidth=1pt,linestyle=dotted}
\pcline(0,4)(6,-2)\Aput{(1,0,1)}
\pcline(6,-2)(0,-8)\Aput{(0,0,1)}
\pcline(0,-8)(-6,-2)\Aput{(1,0,0)}
\pcline(-6,-2)(0,4)\Aput{(0,0,2)}
\psset{linestyle=solid}
\psline(-3,1)(0,0)
\psline[doubleline=true,linewidth=1pt](3,1)(0,0)
\psline(0,0)(0,-4)
\psline[doubleline=true,linewidth=1pt](0,-4)(-3,-5)
\psline(0,-4)(3,-5)
\psarc(0,4){0.3464}{-130}{-50}
}
\rput(80,0){
\psset{linewidth=1pt,linestyle=dotted}
\pcline(0,4)(6,-2)\Aput{(0,0,2)}
\pcline(6,-2)(0,-8)\Aput{(1,0,0)}
\pcline(0,-8)(-6,-2)\Aput{(0,0,2)}
\pcline(-6,-2)(0,4)\Aput{(1,0,0)}
\psset{linestyle=solid}
\psline[doubleline=true,linewidth=1pt](-2.5,1.5)(-1,0)
\psline(2,2)(-1,0)
\psline(-1,0)(-1,-1)
\psline(-1,-1)(-4,-4)
\psline[doubleline=true,linewidth=1pt](-1,-1)(1,-3)
\psline(1,-3)(4,0)
\psline(1,-3)(1,-4)
\psline(1,-4)(-2,-6)
\psline[doubleline=true,linewidth=1pt](1,-4)(2.5,-5.5)
}
\end{pspicture}
\end{center}
This leaves the following four cases
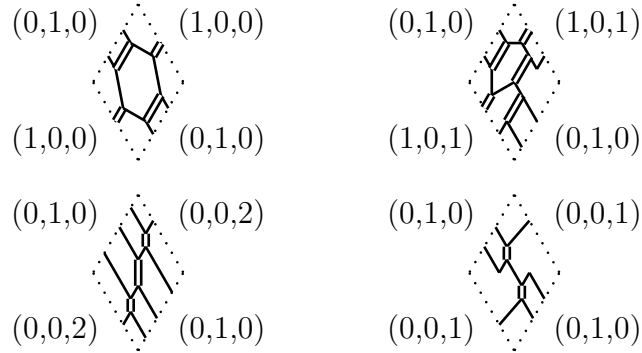
\begin{center}
\psset{xunit=0.1cm,yunit=0.1732cm}
\begin{pspicture}(-6,-8)(56,4)
\rput(0,0){
\psset{linewidth=1pt,linestyle=dotted}
\pcline(0,4)(6,-2)\Aput{(1,0,0)}
\pcline(6,-2)(0,-8)\Aput{(0,1,0)}
\pcline(0,-8)(-6,-2)\Aput{(1,0,0)}
\pcline(-6,-2)(0,4)\Aput{(0,1,0)}
\psset{linestyle=solid}
\psline[doubleline=true,linewidth=1pt](3,1)(2,0)
\psline(4,-4)(3,-3)
\psline(2,-6)(1,-5)
\psline[doubleline=true,linewidth=1pt](-3,-5)(-2,-4)
\psline(-4,0)(-3,-1)
\psline(-2,2)(-1,1)
\psline(2,0)(3,-3)
\psline[doubleline=true,linewidth=1pt](3,-3)(1,-5)
\psline(1,-5)(-2,-4)
\psline(-2,-4)(-3,-1)
\psline[doubleline=true,linewidth=1pt](-3,-1)(-1,1)
\psline(-1,1)(2,0)
}
\rput(50,0){
\psset{linewidth=1pt,linestyle=dotted}
\pcline(0,4)(6,-2)\Aput{(1,0,1)}
\pcline(6,-2)(0,-8)\Aput{(0,1,0)}
\pcline(0,-8)(-6,-2)\Aput{(1,0,1)}
\pcline(-6,-2)(0,4)\Aput{(0,1,0)}
\psset{linestyle=solid}
\psline[doubleline=true,linewidth=1pt](2,2)(1,1)
\psline(3,-1)(2,0)
\psline(1,-3)(0,-2)
\psline[doubleline=true,linewidth=1pt](-4,-4)(-3,-3)
\psline(-4,0)(-3,-1)
\psline(-2,2)(-1,1)
\psline(1,1)(2,0)
\psline[doubleline=true,linewidth=1pt](2,0)(0,-2)
\psline(0,-2)(-3,-3)
\psline(-3,-3)(-3,-1)
\psline[doubleline=true,linewidth=1pt](-3,-1)(-1,1)
\psline(-1,1)(1,1)
\psline(3,-1)(4,0)
\psline[doubleline=true,linewidth=1pt](1,-3)(-1,-5)
\psline(1,-3)(3,-5)
\psline(-1,-5)(1,-7)
\psline(-1,-5)(-2,-6)
}
\end{pspicture}
\end{center}
\begin{center}
\psset{xunit=0.1cm,yunit=0.1732cm}
\begin{pspicture}(-6,-8)(56,4)
\rput(0,0){
\psset{linewidth=1pt,linestyle=dotted}
\pcline(0,4)(6,-2)\Aput{(0,0,2)}
\pcline(6,-2)(0,-8)\Aput{(0,1,0)}
\pcline(0,-8)(-6,-2)\Aput{(0,0,2)}
\pcline(-6,-2)(0,4)\Aput{(0,1,0)}
\psset{linestyle=solid}
\psline[doubleline=true,linewidth=1pt](1,1)(1,0)
\psline(1,0)(0,-1)
\psline[doubleline=true,linewidth=1pt](0,-1)(0,-3)
\psline(0,-3)(-1,-4)
\psline[doubleline=true,linewidth=1pt](-1,-4)(-1,-5)
\psline(1,1)(2,2)
\psline(-1,-5)(-2,-6)
\psline(1,1)(-1,3)
\psline(-1,-5)(1,-7)
\psline(0,-1)(-2.5,1.5)
\psline(0,-3)(2.5,-5.5)
\psline(-1,-4)(-4.5,-0.5)
\psline(1,0)(4.5,-3.5)
}
\rput(50,0){
\psset{linewidth=1pt,linestyle=dotted}
\pcline(0,4)(6,-2)\Aput{(0,0,1)}
\pcline(6,-2)(0,-8)\Aput{(0,1,0)}
\pcline(0,-8)(-6,-2)\Aput{(0,0,1)}
\pcline(-6,-2)(0,4)\Aput{(0,1,0)}
\psset{linestyle=solid}
\psline(-2.5,1.5)(-1,0)
\psline(2,2)(-1,0)
\psline[doubleline=true,linewidth=1pt](-1,0)(-1,-1)
\psline(-1,-1)(-2,-2)
\psline(-2,-2)(-4,0)
\psline(-1,-1)(1,-3)
\psline(1,-3)(2,-2)
\psline(2,-2)(4,-4)
\psline[doubleline=true,linewidth=1pt](1,-3)(1,-4)
\psline(1,-4)(-2,-6)
\psline(1,-4)(2.5,-5.5)
}
\end{pspicture}
\end{center}

Next we associate an element of $\otimes^n\Delta$ to each
of these triangular diagrams with $n$ points on the top edge.
This follows \cite{MR1446615}. First a triangular diagram
interpreted as a tensor is an intertwiner
\[ V(H)\otimes V(D) \rightarrow \otimes^n\Delta \]
so applying this map to the tensor product of highest weight
vectors gives an element of $\otimes^n\Delta$.

\begin{conj}\label{tb}
If we write the tensor associated to a
triangular diagram in terms of the tensor product basis
we get an expression of the form
\begin{equation}\label{tri}
b_D = \sum_{b^\prime \le b} \alpha^{b^\prime}_b b^\prime_T
\end{equation}
where $\alpha^{b^\prime}_b\in \bZ[q,q^{-1}]$ and
$\alpha^{b^\prime}_b=1$ if $b^\prime=b$.
\end{conj}

This means that the tensors associated to the triangular diagrams
are a basis of $\otimes^n\Delta$ and the change of basis matrices
are integral and unitriangular.
The following is a sketch of a proof of this result.

Consider the diagrams in which only some of
the diamonds have been filled in. If we start at the top left
corner of one of these diagrams we have a sequence of alternating
South-East and North-East boundary edges which take us to the
top right corner. Each of these edges is labelled by a dominant
weight. Then the diagram can be considered as map from the
tensor product of the corresponding sequence of highest weight
representations to $\otimes^n\Delta$. Applying this map to the
tensor product of the highest weight vectors gives an element
of $\otimes^n\Delta$. Then we prove that each of these vectors
has the property in \eqref{tri}. The proof is by induction on
the number of diamonds that have been filled in. If no diamonds
have been filled in then this construction just gives the tensor
product of the basis vectors. It is this property that determined
the diagrams we drew in the first place. The inductive step
follows from the observation that if we add one diamond then we 
multiply by a triangular matrix.

A corollary of this conjecture is that the diagrams with $H=0$ and
$D=0$ are a basis of the space of invariant tensors. This corollary
follows from the results in \S \ref{conc} and is part of
Proposition \ref{cp}.

\section{Comparison}\label{conc}
The aim of this section is to explain the relationships between
the constructions of the previous three sections.

Let $A_D(n)$ be the endomorphism algebra of $n$ in $\mathcal{D}$.
There are algebra homomorphisms $A_P(n)\rightarrow A_D(n)$
for $n\ge 0$. This map is defined on the generators by \eqref{hom}
and the defining relations for $A_P(n)$ are satisfied by construction.

In \S \ref{diagrams} we used trivalent graphs to construct the
free strict pivotal category $\widetilde{D}$. Here we take a
more algebraic approach. The full subcategory of $\widetilde{\mathcal{D}}$
whose objects are finite sequences of single edges is generated
as a category by morphisms $A_i\colon n\rightarrow n-2$,
$V_i\colon n-2\rightarrow n$, $K_i\colon n\rightarrow n$
and $H_i\colon n\rightarrow n$. The diagrams for these morphisms
are given in Figure \ref{gens} and
the diagram with index $i$ is obtained by putting
$i-1$ vertical lines on the left and $n-i-1$ vertical lines
on the right.
\begin{figure}
\begin{equation}\label{mon}
\psset{xunit=0.5,yunit=0.5}
\begin{array}{cccc}
A & V & K & H \\
\quad
\begin{pspicture}(2,-1)(4,1)
\psarc(3,-1){0.3}{0}{180}
\end{pspicture}\quad &
\quad
\begin{pspicture}(2,-1)(4,1)
\psarc(3,1){0.3}{180}{360}
\end{pspicture}\quad &
\quad
\begin{pspicture}(2,-1)(4,1)
\psline(2,1)(3,0.5)
\psline(2,-1)(3,-0.5)
\psline[doubleline=true](3,0.5)(3,-0.5)
\psline(3,0.5)(4,1)
\psline(3,-0.5)(4,-1)
\end{pspicture}\quad &
\quad
\begin{pspicture}(2,-1)(4,1)
\psline(2,1)(2.5,0)
\psline(2,-1)(2.5,0)
\psline[doubleline=true](2.5,0)(3.5,0)
\psline(3.5,0)(4,1)
\psline(3.5,0)(4,-1)
\end{pspicture}\quad
\end{array}
\psset{xunit=2,yunit=2}
\end{equation}
\caption{Generators}\label{gens}
\end{figure}

Each morphism is then a composable sequence of these generators and
we define the width of the morphism to be the maximum value of $n$.

\begin{lem}\label{surpt}
A diagram in $A_D(n)$ of width $n$ can be written as a word
in the generators of $A_P(n)$.
\end{lem}
\begin{pf} Since $U_i=V_iA_i$, this means that we can find a 
composable sequence of these generators such that every $V_i$ is
succeeded by $A_i$ and whose diagram is isotopic to the given diagram.
First replace each $V_{n-k}$ by
\[ (V_{n-k}A_{n-k})(V_{n-k+1}A_{n-k+1})\ldots (V_{n-1}A_{n-1})V_n \]
and replace each $A_{n-k}$ by
\[ A_{n}(V_{n-1}A_{n-1})(V_{n-2}A_{n-2})\ldots (V_{n-k}A_{n-k}) \]
Then we draw the diagram in the strip $0\le x\le n+1$ so that for
$1\le m\le n$ the object $m$ the $x$-coordinates of the $m$
points are $1,2,\ldots ,m$. Then in this diagram there is a 
maximum below each minimum. This means that the maxima and minima
are paired by vertical lines which do not intersect the diagram.
Then for each pair we perform an isotopy which shrinks this
vertical line. If each of these vertical lines is sufficiently
short then the resulting diagram has the required form.
\end{pf}

\begin{lem}\label{surpd}
For $n\ge 0$, $A_P(n)\rightarrow A_D(n)$ is surjective.
\end{lem}
\begin{pf}
We show that if $k>0$ then a diagram $D\in A_D(n)$ of width
$n+2k$ can be written as linear combination of diagrams of
width $n+2k-2$. Given $D$ add $k$ minima at the top right
and $k$ maxima at the bottom right. This means we repeat the following
operation until we get a diagram whose width is at most the number
of endpoints on the top and bottom edge
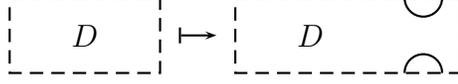
\begin{center}
\begin{pspicture}(0,0)(6,1)
\pspolygon[linestyle=dashed](0,0)(2,0)(2,1)(0,1)
\pspolygon[linestyle=dashed](3,0)(6,0)(6,1)(3,1)
\psarc(5.5,1){0.25}{180}{360}
\psarc(5.5,0){0.25}{0}{180}
\pcline[arrows=|->](2.25,0.5)(2.75,0.5)
\rput(1,0.5){$D$}
\rput(4,0.5){$D$}
\end{pspicture}
\end{center}
This gives a diagram in
$A_D(n+2k)$ of width $n+2k$. This can then be written as a word
in the generators of $A_P(n+2k)$. This word is of the form
$U_{n+2k-1}aU_{n+2k-1}$ for $a\in A_P(n+2k)$. This word is an element
of $A_P(n+2k-2)U_{n+2k-1}$ by Lemma \ref{l} which means that
$D$ has been written as a linear combination of diagrams of width
$n+2k-2$.
\end{pf}

Let $\mathcal{T}_{\bQ}$ be the $\bQ(q)$-linear category of invariant tensors
for the spin representation, $\Delta$. This has objects $0,1,2,\ldots$
and a morphism $n\rightarrow m$ is a linear map
$\otimes^n\Delta\rightarrow\otimes^m\Delta$ which intertwines
the actions of the Drinfeld-Jimbo quantum group $U_q(B_3)$.
This quantum group is defined to be a Hopf algebra over the field
$\bQ(q)$.

Then $\mathcal{T}_{\bQ}$
is a strict pivotal category and so by the universal property
of $\widetilde{\mathcal{D}}$ there is a unique functor of strict
pivotal categories $\widetilde{\mathcal{D}}\rightarrow \mathcal{T}_{\bQ}$
which maps the trivalent vertex to a morphism in $\mathcal{T}_{\bQ}$
from $\Delta\otimes\Delta$ to the vector representation. Note that
this morphism is unique up to multiplication by a non-zero scalar.
\begin{prop}
The defining relations for $\mathcal{D}$ are satisfied so this
functor factorises through $\mathcal{D}$ to give a functor of
$\bZ[\delta]$-linear strict pivotal categories 
$\mathcal{D}\rightarrow \mathcal{T}_{\bQ}$.
\end{prop}
This result is the motivation for the relations.
\begin{pf} Let $B_3$ be
the three string braid group. Then $B_3$ acts on the invariant 
tensors in $\otimes^3\Delta$. The dimensions of these representations
are known from the decomposition of $\otimes^3\Delta$. Also the
eigenvalues of the standard generators are known. Since each
representation has dimension at most five this information
determines the representations by \cite{MR1815266}. This gives
the representations \S \ref{rep}. Given these representations and
using the pivotal structure one can check that the defining relations
for $\mathcal{D}$ are satisfied.
\end{pf}

Let $A_T(n)_{\bQ}$ be the endomorphism algebra of $n$ in $\mathcal{T}_{\bQ}$.
Then in the remainder of this section we will use the homomorphisms
\begin{equation}\label{homa}
A_P(n)\rightarrow A_D(n)\rightarrow A_T(n)_{\bQ}
\end{equation}
to compare these algebras.

Let $\mathcal{D}_{\bQ}$ be the free $\bQ(q)$-linear category on $\mathcal{D}$
so $\mathcal{D}_{\bQ}=\mathcal{D}\otimes_{\bZ[\delta]}\bQ(q)$.
Then for $n\ge 0$ we can define $A_D(n)_{\bQ}$ either as the endomorphism
algebra of the object $n$ of $\mathcal{D}_{\bQ}$ or as 
$A_D(n)\otimes_{\bZ[\delta]}\bQ(q)$.
Then the functor $\mathcal{D}\rightarrow \mathcal{T}_{\bQ}$ extends
to a functor $\mathcal{D}_{\bQ}\rightarrow \mathcal{T}_{\bQ}$ and induces an
algebra homomorphism $A_D(n)_{\bQ}\rightarrow A_T(n)_{\bQ}$ for each $n\ge 0$.

For $n> 1$ we have generators $U_1,\ldots ,U_{n-1}\in A_P(n)$.
We also consider these as elements of the algebra $A_P(n)_{\bQ}$
and as elements of $A_T(n)_{\bQ}$ using the homomorphisms in
\eqref{homa}. Then we define $H_P(n)_{\bQ}$ to be the quotient of
$A_P(n)_{\bQ}$ by the ideal generated by these elements and we define
$H_T(n)_{\bQ}$ to be the quotient of $A_T(n)_{\bQ}$ by the ideal
generated by these elements. Then, for $n>0$ we have an induced
homomorphism $H_P(n)_{\bQ}\rightarrow H_T(n)_{\bQ}$.

In \cite{MR1090432} it is shown that these homomorphisms are
isomorphisms. Furthermore these algebras are direct sums of matrix
algebras and the Bratteli diagram is the same as the Bratteli diagram
of the centraliser algebras of the tensor powers of the four dimensional
fundamental representation of the simple Lie algebra of type $B_2$.
This representation can be taken to be the vector representation
of $\fsp(4)$ or the spin representation of $\fso(5)$.

\begin{prop}\label{isopt} For $n>0$, the homomorphism
$A_P(n)_{\bQ}\rightarrow A_T(n)_{\bQ}$ is surjective.
\end{prop}
\begin{pf} The proof is by induction on $n$. We have
\[ A_T(n+1)_{\bQ} \cong A_T(n)_{\bQ}U_nA_T(n)_{\bQ}\oplus H_T(n)_{\bQ} \]
The inductive step now follows from the inductive hypothesis and the
result that the homomorphism $H_P(n)_{\bQ}\rightarrow H_T(n)_{\bQ}$
is surjective.
\end{pf}

\begin{cor}\label{ptiso} For $n>0$, the homomorphism
$A_P(n)_{\bQ}\rightarrow A_T(n)_{\bQ}$ is an isomorphism.
\end{cor}
\begin{pf} The number of irreducible words in Definition \ref{irr}
is also $\dim A_T(n)_{\bQ}$. Hence the observation that the
set of irreducible words is a spanning set for $\dim A_P(n)_{\bQ}$
shows that $\dim A_P(n)_{\bQ}\le \dim A_T(n)_{\bQ}$.
\end{pf}

\begin{cor}\label{iso} For all $n>0$, the homomorphism
$A_P(n) \rightarrow A_D(n)$ is an isomorphism.
\end{cor}
\begin{pf} Corollary \ref{ptiso} implies that this homomorphism
is injective. The result follows from this observation and
Lemma \ref{surpd}.
\end{pf}

\begin{cor} The strict pivotal functor
$\mathcal{D}_{\bQ}\rightarrow\mathcal{T}_{\bQ}$ is an isomorphism.
\end{cor}
\begin{pf} The categories $\mathcal{D}$ and $\mathcal{T}$ are strict
pivotal and so if $n+m=2p$ we have canonical isomorphisms
$\Hom_\mathcal{D}(n,m)\rightarrow A_D(p)$ and
$\Hom_\mathcal{T}(n,m)\rightarrow A_T(p)$. Since the functor
is strict pivotal and induces an isomorphism $A_D(p)\rightarrow A_T(p)$
it also induces an isomorphism
$\Hom_\mathcal{D}(n,m)\rightarrow \Hom_\mathcal{T}(n,m)$.

If $n+m$ is odd then $\Hom_\mathcal{D}(n,m)$ and $\Hom_\mathcal{T}(n,m)$
are both the zero module so there is nothing to prove.
\end{pf}

This result can be rephrased to say that the pivotal category $\mathcal{D}$
is an integral form or order for the category $\mathcal{T}_{\bQ}$.
Another construction of an integral form for $\mathcal{T}_{\bQ}$ is given
in \cite[Part IV]{MR1227098}. Then Conjecture \ref{tb} implies that these two
integral forms are equivalent.

\section{Confluence}\label{conf}
In this section we apply the theory of rewrite rules to the presentations
of the algebras $A_P(n)$. This theory originates from the diamond lemma
in \cite{MR0007372} and our account is based on \cite[Chapter 2]{MR1267733}.
This theory has been applied to the Hecke and Temperley-Lieb algebras
in \cite{MR1925138}. Our approach differs from these standard approaches
in that we take $A_P(n)$ to be a quotient of the free algebra on a
commutation monoid (instead of a free monoid) and we only require a
reduction order to be a partial order (instead of a linear order).

Let $X$ be a set and $W\subset X\times X$ a relation.
We write $x\rightarrow y$ if $(x,y)\in W$.
Let the relation $\stackrel{*}{\rightarrow}$ be
the reflexive transitive closure of the relation $W$.
Let $\simeq$ be the equivalence relation generated by $W$.

The first result is the following.
\begin{defn}\label{new1} A rewrite system is confluent if either
of the following two equivalent conditions is satisfied
\begin{enumerate}
\item If $u\simeq v$ then there exists an $x$ such that
$u\stackrel{*}{\rightarrow}x$ and $v\stackrel{*}{\rightarrow}x$.
\item If $u\stackrel{*}{\rightarrow}x$
and $u\stackrel{*}{\rightarrow}y$ then there exists a $v$
such that $x\stackrel{*}{\rightarrow}v$ and $y\stackrel{*}{\rightarrow}v$.
\end{enumerate}
A rewrite system is locally confluent if whenever
$u{\rightarrow}x$ and $u{\rightarrow}y$ then there exists a $v$
such that $x\stackrel{*}{\rightarrow}v$ and $y\stackrel{*}{\rightarrow}v$.
\end{defn}

It is clear that a confluent rewrite system is locally confluent
but there are examples of locally confluent rewrite systems that
are not confluent.
\begin{defn} A rewrite system is terminal if there is no infinite
sequence $x_0,x_1,\ldots $ such that $x_{i-1}\rightarrow x_i$
for all $i>0$.
\end{defn}
Note that if a rewrite system is terminal then the relation
$\stackrel{*}{\rightarrow}$ is a partial order.
\begin{prop}\label{new2} A rewrite system that is terminal and
locally confluent is confluent.
\end{prop}

Now we apply this to finitely presented algebras.
Let $K$ be a commutative ring and $M$ a monoid. 
Let $KM$ be the monoid algebra of $M$ over $K$. Then every element
$u\in KM$ can be written uniquely as $\sum_{m\in S} r_mm$
where $S$ is a finite subset of $M$ and $r_m\in K\setminus \{0\}$.
The subset $S$ is called the support of $u$ and is denoted by
$\supp(u)$.

Let $R$ be a finite set of ordered pairs $(p,P)$ where $p\in M$,
$P\in KM$ and $p\notin\supp(P)$. Then let $I\subset KM$ be the ideal
generated by
the set $\{p-P|(p,P)\in R\}$ and let $A$ be the $K$-algebra
$KM/I$.

Now define a relation $W$ on the set $KM$ by
\[ rxpy+u \rightarrow rxPy+u \]
where $r\in K\setminus \{0\}$, $x,y\in M$, $u\in KM$, $(p,P)\in R$
and $p\notin \supp(u)$. Then the equivalence relation $u\simeq v$
is the equivalence relation $u-v\in I$. 

\begin{defn}\label{irp} An element $u\in M$ is reducible if it can be
written as $u=apb$ where $a,b\in M$ and for some $P\in KM$
$(p,P)\in R$. An element $u\in M$ is irreducible
if it is not reducible.
\end{defn}

\begin{lem} The rewrite system $W$ is confluent if and only if
$A$ is the free $K$-module on the set of irreducible elements
of $M$.
\end{lem}

The following Proposition \ref{ov} gives a criterion which can
be checked by a finite calculation and which implies that the
rewrite system is locally confluent.

\begin{defn} An overlap consists of elements $u,v,w\in M$
together with elements $P,Q\in KM$ such that $v\ne 1$,
$(uv,P)\in R$ and $(vw,Q)\in R$. This overlap is unambiguous
if there is an $x\in KM$ such that $Pw\stackrel{*}{\rightarrow}x$ and
$uQ\stackrel{*}{\rightarrow}x$.
\end{defn}

\begin{prop}\label{ov} The rewrite system $W$ is locally confluent
if and only if every overlap is unambiguous.
\end{prop}

\begin{defn} A reduction order on $M$ is a partial order which
is invariant under both left and right translations and such that,
for each $m\in M$, the set $\{m^\prime\in M|m^\prime <m\}$ is finite.
\end{defn}

Then the standard method of showing that a rewrite system is terminal
is to construct a reduction order on $M$ such that for each
$(p,P)\in R$ we have $a<p$ for each $a\in\supp(P)$.

Next we apply this theory to the presentations of the algebras
$A_P(n)$.
\begin{defn}\label{com}
For $n\ge 1$ let $C(n)$ be the monoid with the same set of generators
as $A_P(n)$ and with defining relations the commuting relations.
That is, for $1\le i,j\le n-1$ with $|i-j|>1$,
if $a\in\{U_i,K_i,H_i\}$ and  $b\in\{U_j,K_j,H_j\}$ then $ab=ba$.
\end{defn}

Let the type of a word in the generators be the sequence of subscripts.
Construct a rewrite system for $A_P(3)$ by taking the set of pairs
$(p,P)$ where the words $p$ are the words of type $(1,1)$ and $(2,2)$,
the words of type $(2,1,2)$ and the words of type $(1,2,1)$ which are
not irreducible words. Then, for each $p$, we can use the relations
to write $p$ uniquely as a linear combination of irreducible words.
Then, by construction, this is a confluent rewrite system for $A_P(3)$.

Then, for $n>3$, we obtain a rewrite system by taking the union over
$1\le i\le n-2$ of the sets obtained by changing each subscript $1$
to $i$ and each subscript $2$ to $i+1$. This gives a finite set of
rewrite rules, $R(n)$. Although we have described these rewrite rules
using words, all of these words are well-defined elements of $C(n)$

Define a partial order on $C(n)$ recursively. This definition is similar
to the definition of a wreath product order on a free monoid. On the free monoid
generated by $\{U_i,K_i,H_i\}$ take the length plus reverse lexicographic
order with the generators ordered by $U_i<K_i<H_i$. Let $w$ be a word
in the generators of $C(n+1)$. Then we obtain a word $u$ in the generators
$\{U_n,K_n,H_n\}$ by deleting all generators whose subscript is not $n$
and a word $v$ in the generators of $C(n)$ by deleting all generators whose
subscript is $n$. Now note that $u$ and $v\in C(n)$ depend only on
$w\in C(n+1)$ and not on the choice of word representing $w$.
Then let $w_1,w_2\in C(n+1)$ with corresponding words $u_1,u_2$ and
elements $v_1,v_2\in C(n)$. Then we define $w_1<w_2$ if $u_1<u_2$ or
$u-1=u_2$ and $v_1<v_2$.

Then this is a reduction order and has the property that if $(p,P)\in R(n)$
and $a\in\supp(P)$ then $a<p$. This shows that the rewrite system $W(n)$
associated to $R(n)$ is terminal. Note that the set of irreducible
words defined in Definition \ref{irr} coincides with the set of
irreducible elements for the rewrite system $W(n)$ in Definition \ref{irp}.
It now follows from the theory of rewrite systems that the
following are equivalent, for any $n>3$:
\begin{enumerate}
\item The rewrite system $W(n)$ for $A_P(n)$ is confluent.
\item The rewrite system $W(n)$ for $A_P(n)$ is locally confluent.
\item The $\bZ[\delta]$-module underlying $A_P(n)$ is the
free module on the set of irreducible words.
\end{enumerate}
The third property has been proved in Proposition \ref{isopt}
using the representation theory of the quantum group $U_q(B_3)$
and so we conclude that these rewrite systems are confluent.

However it is also possible to prove that these rewrite systems
are locally confluent by directly checking that all overlaps
are unambiguous. This depends on the observation that this
holds for any $n\ge 4$ if and only if it holds for $n=4$.
This holds since any overlap can only involve subscripts of the
form $i-1$, $i$ and $i+1$ for some $i>1$. This case can in principle
be checked directly since it is a finite calculation. This case
can also be checked indirectly by giving a finite dimensional
representation in which the irreducible words are linearly independent.
This case can also be checked indirectly using a computer algebra
package.

\section{Cellular algebras}\label{cell}
First we extend the definition of a cellular algebra given in
\cite{MR1376244} to cellular categories. The basic example
is the Temperley-Lieb category discussed in \cite{MR1343661}.
The proof we give here also applies to the categories of diagrams
constructed from the rank two simple Lie algebras
in \cite{MR1403861}. Further examples are the affine Temperley-Lieb
category in \cite{MR1659204} and the partition category whose
endomorphism algebras are the partition algebras.

\begin{defn}\label{cc}
Let $R$ be a commutative ring with identity. Let $\cA$ be a $R$-linear
category with an anti-involution $*$. Then cell datum for $\cA$ consists
of a partially ordered set $\Lambda$, a finite set $M(n,\lambda)$ 
for each $\lambda\in\Lambda$ and each object $n$ of $\cA$, and for 
$\lambda\in\Lambda$ and $n$,$m$ any two objects of $\cA$ we have an inclusion
\[ C\colon M(n,\lambda)\times M(m,\lambda)\rightarrow \Hom_\cA(n,m) \]
\[ C\colon (S,T)\mapsto C^\lambda_{S,T} \]
The conditions that this datum is required to satisfy are:
\begin{enumerate}
\item[C-1] For all objects $n$ and $m$, the image of the map
\[ C\colon \coprod_{\lambda\in\Lambda} M(n,\lambda)\times M(m,\lambda) 
\rightarrow \Hom_\cA(n,m) \]
is a basis for $\Hom_\cA(n,m)$ as an $R$-module.
\item[C-2] For all objects $n,m$, all $\lambda\in\Lambda$ and all
$S\in M(n,\lambda)$, $T\in M(m,\lambda)$ we have
\[ \left(C^\lambda_{S,T}\right)^*=C^{\lambda}_{T,S} \]
\item[C-3] For all objects $p,n,m$, all $\lambda\in\Lambda$ and all
$a\in \Hom_\cA(p,n)$,
$S\in M(n,\lambda)$, $T\in M(m,\lambda)$ we have
\[ aC^\lambda(S,T)=\sum_{S^\prime\in M(p,\lambda)}
r_a(S^\prime,S)C^\lambda_{S^\prime ,T}
\mod \cA(<\lambda) \]
where $r_a(S^\prime,S)\in R$ is independent of $T$ and 
$\cA(<\lambda)$ is the $R$-linear span of
\[ \left\{ C^\mu_{S,T} | \mu < \lambda ; S\in M(p,\mu),T\in M(m,\mu)
\right\} \]
\end{enumerate}
\end{defn}
Another way of formulating C-3 is that for all objects
$p,n,m$, all $\lambda\in\Lambda$ and all
$S\in M(p,\lambda)$,$T\in M(n,\lambda)$,
$U\in M(n,\lambda)$,$V\in M(m,\lambda)$ we have
\[ C^\lambda(S,T)C^\lambda(U,V)=\langle T,U\rangle C^\lambda_{S,V}
\mod \cA(<\lambda) \]
where $\langle T,U\rangle$ is independent of $S$ and $V$.

A consequence of this definition is that for all $\lambda\in\Lambda$
we have ideals $\cA(<\lambda)$ and $\cA(\le\lambda)$ where
for all objects $n,m$ the subspace $\cA(<\lambda)(n,m)$ is the
$R$-linear span of 
\[ \left\{ C^\mu_{S,T} | \mu < \lambda ; S\in M(n,\mu),T\in M(m,\mu) \right\} \]
and the subspace $\cA(\le\lambda)(n,m)$ is the $R$-linear span of
\[ \left\{ C^\mu_{S,T} | \mu \le \lambda ; S\in M(n,\mu),T\in M(m,\mu) \right\} \]

This is a generalisation of the definition of a cellular algebra
since we can regard any algebra over $R$ as an $R$-linear category
with one object.

More significantly, if $\cA$ is a cellular category then $\End(n)$
is a cellular algebra for any object $n$ of $\cA$.

A functor $\phi\colon\cA\rightarrow\cA^\prime$ between cellular categories
is cellular if we have a map of partially ordered sets
$\phi\colon \Lambda\rightarrow\Lambda^\prime$ and set maps
\[ \phi_\lambda\colon M(n,\lambda)\rightarrow M^\prime(\phi(n),\phi(\lambda)) \]
such that
\[ \phi( C^\lambda_{S,T})= C^{\phi(\lambda)}_{\phi_\lambda(S),\phi_\lambda(T)}
\mod \cA^\prime(<\phi(\lambda)) \]for all $\lambda$ and all $S,T$.

One reason for considering cellular categories instead of just 
the cellular endomorphism algebras is that for each object $n$
and each $\lambda\in\Lambda$ we have an $R$-linear functor
$\rho(n;\lambda)$ from $\cA$ to the category of left
$\End_\cA(n)$-modules which on objects is given by
\[ \rho(n;\lambda)\colon m\mapsto \Hom_{\cA(\lambda)}(n,m) \]
Furthermore if $\mu<\lambda$ then we have a natural transformation
from $\rho(n;\mu)$ to $\rho(n;\lambda)$. This is used
in \cite{MR1659204}.

\begin{prop}\label{cp} The category $\mathcal{D}$ with the $*$-functor given
by the pivotal structure is cellular.
\end{prop}

\begin{pf} The partially ordered set $\Lambda$ is the set of dominant
weights. For $n$ an object of $\mathcal{D}$ and $\lambda$ a dominant
weight we take the set $M(n,\lambda)$ to be a set of diagrams
(or morphisms in $\widetilde{\mathcal{D}}$). We require that each diagram
in $M(n,\lambda)$ has $n$ boundary points on the top edge and that the
bottom edge is a minimal cut path of weight $\lambda$
(see Definition \ref{cut}). Furthermore we require that the bottom
edge is minimal with these properties. This means that each diagram
$D\in M(n,\lambda)$ has the property that if $D=D_1D_2$ is a factorisation
in $\widetilde{\mathcal{D}}$ such that the bottom edge of $D_1$ has
weight $\lambda$ then $D=D_1$.

Next we construct the map $C$. Let $D^\prime\in M(n,\lambda)$ and
$D\in M(m,\lambda)$. Then we define $C(D^\prime ,D)$ to be
the composition in $\widetilde{\mathcal{D}}$, $D^\prime HD^*$
where the diagram $H$ is uniquely determined by the properties
that this composition exists and gives an irreducible diagram such
that every cut path that traverses the rectangle has weight at least
$\lambda$.

This gives the data in Definition \ref{cc}. Next we verify that the
conditions are satisfied. Each irreducible $(n,m)$-diagram can be
written uniquely as $C(D^\prime ,D)=D^\prime HD^*$. The weight
$\lambda$ is the minimum weight of a cut path which traverses the
rectangle. This shows that the diagram has a factorisation as
$D^\prime D^*$ where the bottom edges of $D$ and $D^\prime$ are
minimal cut paths of weight $\lambda$. Then to obtain the
factorisation as $C(D^\prime,D)=D^\prime HD^*$ apply the analogue
of \cite[Lemma 6.6]{MR1403861} to both the top and bottom edge of the
rectangle.

It has been shown in \S \ref{conc} that the set of
irreducible $(n,m)$-diagrams is a basis of $\Hom_{\mathcal{D}}(n,m)$
and hence property C-1 holds.

The property C-2 holds by construction since
\[ C(D,D^\prime)^*=(D^\prime HD^*)^*=DH^*(D^\prime)^*=C(D^\prime,D) \]
The property C-3 holds by inspecting the relations for $\mathcal{D}$.
\end{pf}
\begin{cor} For $n>0$ the algebra $A_P(n)\cong A_D(n)$ is cellular
over $\bZ[\delta]$.
\end{cor}

Let $\Bbbk$ be a field and $\phi\colon \bZ[\delta]\rightarrow \Bbbk$
a ring homomorphism. Then by applying the functor
$-\otimes_{\bZ[\delta]}\Bbbk$ to the algebra $A_P(n)\cong A_D(n)$
we obtain a $\Bbbk$-algebra which we denote by $A(n)_\phi$.
In the rest of this section we will write $[r]$ for $\phi([r])\in\Bbbk$.

The definition of a quasi-hereditary algebra is given in
\cite{MR961165} and cellular algebras over a field which are
quasi-hereditary are characterised in \cite{MR1721800}.
\begin{prop} If $\delta(\delta^2-2)(\delta^2-3)(\delta^4-5\delta^2+5)\ne 0$
then $A(n)_\phi$ has the following properties
\begin{enumerate}
\item The algebra $A(n)_\phi$ is quasi-hereditary.
\item The algebra $A(n)_\phi$ has finite global dimension.
\item The Cartan matrix of $A(n)_\phi$ has determinant one.
\end{enumerate}
\end{prop}
\begin{pf} We apply \cite[Lemma 2.1 (3)]{MR1721800} to show that
the algebra $A(n)_\phi$ is quasi-hereditary. The result is then an
application of \cite[Theorem 3.1]{MR1721800}.

Let $0=J_0\subset J_1\subset \ldots \subset J_N=A(n)_\phi$ be a
cell chain of ideals. Then we show that for $1\le r\le N$ that
there is an idempotent $e$ such that $e\in J_r$ and $e\notin J_{r-1}$.
It follows that $J_r^2\nsubseteq J_{r-1}$.

For $1\le i\le n-1$ we have orthogonal idempotents
\[ \frac{1}{[2]^2}U_i, \frac{[3]}{[2]^2[6]}K_i,
\frac{1}{[2]^2}+\frac{1}{[2]^2}H_i+\frac{1}{[2]^3}K_i
-\frac{[4][5]}{[2]^3[10]}U_i \]
which we denote by $u_i$,$k_i$ and $e_i$. Then for each dominant
weight $(a,b,c)$ such that $2a+2b+2c\le n$ we have an idempotent
\[ (u_1u_3\ldots u_{2a-1})(k_{2a+1}k_{2a+3}\ldots k_{2a+2b-1})
(e_{2a+2b+1}e_{2a+2b+3}\ldots e_{2a+2b+2c-1}) \]
\end{pf}
%\bibliographystyle{elsart-harv}
%\nocite{*}
%\bibliography{spin}

\end{document}